\documentclass{amsart}
\usepackage{geometry}
\usepackage{amssymb}
\usepackage{amsfonts}
\usepackage{amsmath}
\usepackage{graphicx}
\usepackage{color}
\usepackage{cite}
\usepackage{chngcntr}

\counterwithout{figure}{section}

\theoremstyle{definition}
\newtheorem{thm}{Theorem}[section]
\newtheorem{example}[thm]{Example}

\newcommand{\cL}{\mathcal{L}}
\newcommand{\cP}{\mathcal{P}}
\newcommand{\bLambda}{\bar{\Lambda}}
\newcommand{\bPi}{\bar{\Pi}}
\newcommand{\bcL}{\bar{\cL}}
\newcommand{\bcP}{\bar{\cP}}

\newcommand{\ffac}[2]{\left( #1 \right)_{\underline{#2}}}
\newcommand{\rfac}[2]{\left( #1 \right)_{\overline{#2}}}

\begin{document}

\title[Equivalence relations for Laguerre and Jacobi pseudowronskians]{Shape invariance and equivalence relations for pseudowronskians of Laguerre and Jacobi polynomials}

\author{David G\'omez-Ullate}
\address{Instituto de Ciencias Matem\'aticas (CSIC-UAM-UC3M-UCM),  C/ Nicolas Cabrera 15, 28049 Madrid, Spain.}
\address{Departamento de F\'isica Te\'orica, Universidad Complutense de
Madrid, 28040 Madrid, Spain.}
\author{ Yves Grandati}
\address{Laboratoire de Physique et Chimie Th\'eoriques, UMR CNRS 7019, Universit\'e de Lorraine--Site de Metz,
1 Bvd D. F. Arago, F-57070, Metz, France.}
\author{Robert Milson}
\address{Department of Mathematics and Statistics, Dalhousie University,
Halifax, NS, B3H 3J5, Canada.}
\email{david.gomez-ullate@icmat.es, grandati@univ-metz.fr, rmilson@dal.ca}

\begin{abstract}
In  a previous paper we derived equivalence relations for pseudo-Wronskian determinants of Hermite polynomials.
In this paper we obtain the analogous result for Laguerre and Jacobi polynomials.
The equivalence formulas are richer in this case since rational Darboux transformations can be defined for four families of seed functions, as opposed to only two families in the Hermite case. The pseudo-Wronskian determinants of Laguerre and Jacobi type will thus depend on two Maya diagrams, while Hermite pseudo-Wronskians depend on just one Maya diagram.
We show that these equivalence relations can be interpreted as the general transcription of shape invariance and specific discrete symmetries acting on the parameters of the isotonic oscillator and Darboux-P\"{o}schl-Teller potential.
\end{abstract}

\keywords{Wronskian determinants, exceptional orthogonal polynomials, Maya diagrams, Durfee symbols, shape invariant potentials, Darboux transformations}
\subjclass[2010]{33C45, 81Q80, 42C05}

\maketitle

%%%%%%%%%%%%%%%%%%%%%%%%%%%%%%%%%%%%%%%%%%%%%%%%%%%%%%%%%%%%%%%%%%%
\section{Introduction}
%%%%%%%%%%%%%%%%%%%%%%%%%%%%%%%%%%%%%%%%%%%%%%%%%%%%%%%%%%%%%%%%%%%

Wronskian determinants arise in a natural fashion when iterating Darboux transformations, \cite{crum}. These ideas have been extensively used in the theory of integrable systems generally known as the \textit{dressing method} (see for instance \cite{Freeman1983, matveev} and references therein).

In the context of classical orthogonal polynomials, Darboux transformations are essentially a factorization of the differential operator, and lead directly to \textit{exceptional orthogonal polynomials}, a new class of Sturm-Liouville polynomial families originally introduced in \cite{gomez,gomez1}, and developed over the past decade by many authors \cite{Quesne2009,Odake2009a, Gomez-Ullate2015b, Odake2011c,Marquette2013a,Duran2014a,Duran2014b,Post2012a,Gomez-Ullate2012a, Gomez-Ullate2013a,Gomez-Ullate2012,Grandati2011b,Kuijlaars2014,Haese-Hill}.

From a physical point of view, exceptional orthogonal polynomials are (up to a gauge factor) the eigenfunctions of exactly solvable quantum potentials obtained by applying rational Darboux transformations  to the harmonic oscillator (Hermite), isotonic oscillator (Laguerre) and trigonomeric Darboux-P\"{o}schl-Teller potential (Jacobi). The new exactly solvable potentials are known are rational extensions \cite{GGM,grandati5,GCB}, since the original potentials are modified by the addition of rational terms in a suitable variable.

In the mathematical context, Wronskian determinants whose entries are classical polynomials were considered by Karlin and Szeg\H{o} \cite{Karlin1960}, who provide expressions for the number of real zeros of the Wronskian of a sequence of consecutive polynomials. Sequences for which the Wronskian has constant sign  were characterized by Adler in \cite{Adler1994}, and the result for arbitrary sequences has been completed by Garc\'ia-Ferrero et al. \cite{Garcia-Ferrero2015}. The complex zeros of these Wronskians form regular patterns and have also been the subject of intense study, \cite{clarkson2009,Felder2012}

These Wronskian determinants of classical polynomials play an important role in the construction of rational solutions to nonlinear integrable equations of Painlev\'e type \cite{veselov93,FP2008,kajiwara,bermudez}.  All known rational solutions to PIV, PV, and their higher order generalizations known as $A_N$-Painlev\'e or  Noumi-Yamada equations seem to be expressible in terms of suitable Wronskians whose entries are Hermite and Laguerre polynomials, although a rigorous proof of this fact has not been given yet (see \cite{noumi} and references therein for a review of the symmetry approach to Painlev\'e equations).

In this paper we are concerned with describing the existence of an infinite number of identities among Wronskian determinants whose entries are Laguerre and Jacobi polynomials. Some of these equivalences were already noted in \cite{odake,GGM1} at the Schr\"odinger's picture, but it was not until \cite{GGM4} that a full description was given for the Hermite case. From the quantum mechanical perspective, the identities reflect a notion of equivalence among different sets of Darboux transformations that end up in the same transformed potential, up to an additive shift.
The term \textit{pseudo-Wronskian} was introduced in \cite{GGM4} to describe the set of rules to construct the determinants with polynomial entries in terms of a Maya diagram, a type of diagram originally introduced by Sato in the theory of integrable systems \cite{ohta} and first applied to exceptional orthogonal polynomials by Takemura, \cite{takemura}. A Hermite pseudo-Wronskian is indexed by a single Maya diagram (or equivalently, a single partition), corresponding to the fact that a rational Darboux transformation of the harmonic oscillator is indexed by a finite sequence of integers specifying the seed functions for the transformation. The equivalance at the level of Maya diagrams corresponds to shifting the origin of the diagram. 

In the Laguerre and Jacobi cases discussed here, the situation is richer as two independent families of seed functions can be considered for rational Darboux transformations, and thus two Maya diagrams must be given to specify the set of seed funtions in each family. Instead of a partition, these rational Darboux transformations are indexed by a Universal Character, which is essentially a pair of partitions as defined by Koike \cite{koike}. Universal characters generalize Schur polynomials, and they have been used in the construction of rational solutions to $A_N$-Painlev\'e equations, \cite{tsuda,tsuda2}.

Although we do not pursue this direction in this paper, it is worth noting that the pseudo-Wronskian determinants discussed here can be regarded as an extension of Jacobi-Trudi formulas in the theory of symmetric functions, \cite{FH}. Original Jacobi-Trudi formulas express a Schur polynomial $s_\lambda$ associated to a given partition $\lambda$ as a determinant whose
entries are complete homogeneous symmetric polynomials. Jacobi-Trudi
formulas were extended to classical orthogonal polynomials in \cite{Sergeev2014} and they are connected to exceptional orthogonal polynomials, \cite{Grandati-Schur}.  The pseudo-Wronskian determinants in this paper can be regarded as a generalization of Jacobi-Trudi formulas that involve not just the partition $\lambda$ or its conjugate partition $\lambda'$, but also intermediate representations described by Durfee symbols, \cite{andrews1}.

The paper is organized as follows. In Section~\ref{sec:PUC} we introduce all the notions in combinatorics needed to formulate our results: Maya diagrams, Durfee symbols, etc. together with some basic operations on them.
In Sections~\ref{sec:DT} and \ref{sec:TSI} we review the necessary concepts on Darboux transformations and shape invariant potentials. It includes a description of the discrete $\Gamma$ symmetries acting on the parameters of the potentials that allow to define the notions of extended and shadow spectra, i.e. the families of possible seed functions for a rational Darboux transformation.
The action of iterated Darboux transformations on the isotonic oscillator and Darboux-P\"{o}schl-Teller potential are described in Section~\ref{sec:shapeinv}, and they are indexed by punctured universal characters.
Finally, we introduce the Laguerre and Jacobi pseudo-Wronskians in Sections~\ref{sec:L} and \ref{sec:J}, together with the equivalence formulas and some explicit examples.

We would like to note that some of the equivalence formulas reported in this paper for the Laguerre pseudo-Wronskians have been obtained independently by Bonneux and Kuijlaars, \cite{Bonneux}, albeit with a different choice of indexing between seed functions and Maya diagrams.

%%%%%%%%%%%%%%%%%%%%%%%%%%%%%%%%%%%%%%%%%%%%%%%%%%%%%%%%%%%%%%%%%%%
\section{Punctured Young diagrams and Durfee rectangles} \label{sec:PUC}
%%%%%%%%%%%%%%%%%%%%%%%%%%%%%%%%%%%%%%%%%%%%%%%%%%%%%%%%%%%%%%%%%%%

\subsection{Maya diagrams}

We define a \textit{Maya diagram} as an infinite row of boxes labelled by
relative integers and which can be empty or filled by at most one "particle"
(graphically represented by a filled dot). All the boxes sufficiently far away to the left are filled and all the boxes sufficiently far away to the
right are empty.

\begin{figure}[ht]
\centering
\includegraphics[width=10cm]{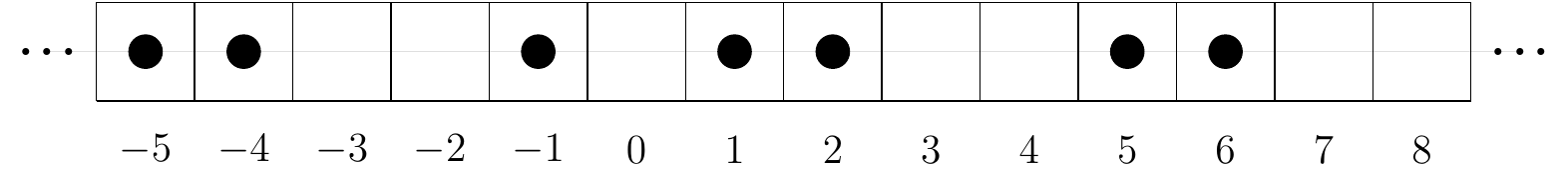}
\caption{Unique Maya diagrama associated to the tuple $(6,5,2,1,-2-3)$. }
\label{fig:maya1}
\end{figure}

A tuple of relative integers, $N_{m}=\left( n_{1},...,n_{k},\overline{n}_{%
\overline{k}},...,\overline{n}_{1}\right) $, with $n_{i}>n_{i+1}\geq 0$ and $%
0>\overline{n}_{i+1}>\overline{n}_{i}$ $(m=k+\overline{k})$, can be
represented in an unique way by a Maya diagram, as depicted in Figure \ref{fig:maya1}. Positive integers in the tuple correspond to filled boxes to the right of the origin and negative integers in the tuple to empty boxes to the left of the origin.

In the following we use the term Maya diagram to designate both the diagram
in graphical form and for the associated tuple. If the set of $%
\overline{n}_{i}$ (respectively $n_{i}$) is empty, then the Maya diagram is
said to be positive (respectively, negative). In the case
of a positive Maya diagram, if all the entries are non-zero, then the Maya
diagram is said to be \textit{strictly positive.}

\subsection{Punctured Young diagrams}

The correspondence between Maya and Young diagrams is well known \cite%
{noumi,hirota}. Progressing from the left to the right, to each filled box
we associate a vertical step upward and to each empty box an horizontal step
to the right. The obtained broken line then serves as the lower and right
boundary of the corresponding Young diagram.

\begin{figure}[h]
\centering
\begin{tabular}{cc}
\includegraphics[width=0.45\textwidth]{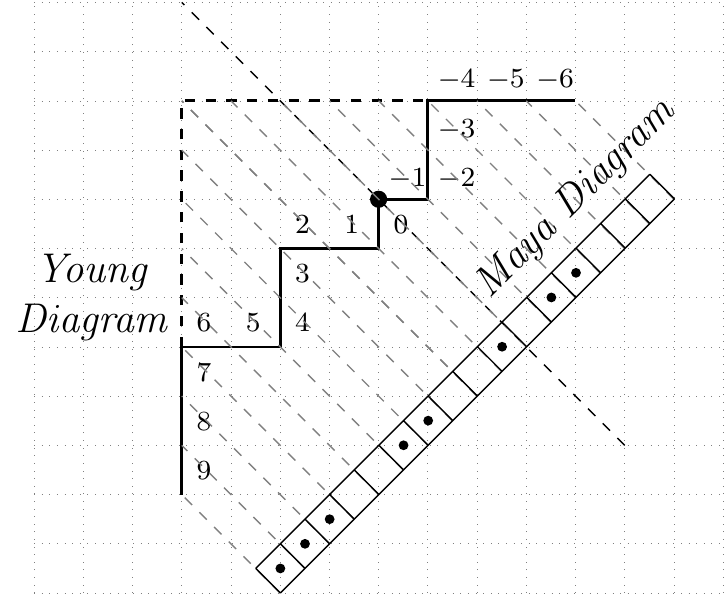} & 
\includegraphics[width=0.45\textwidth]{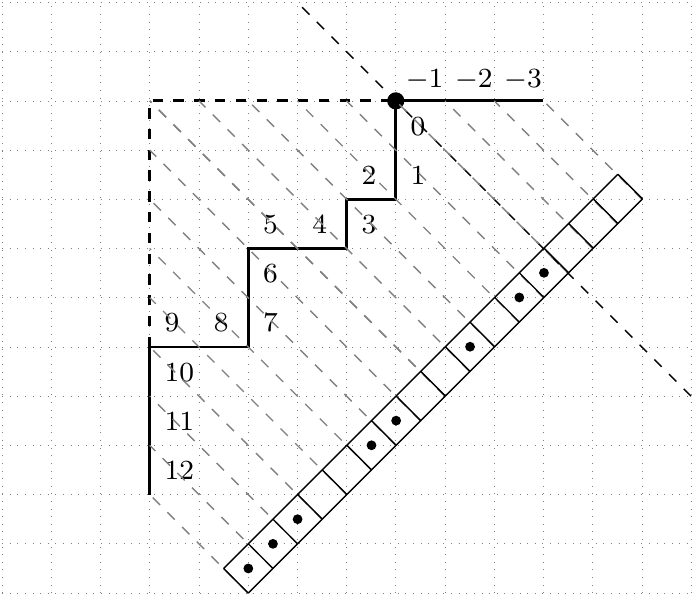} 
\end{tabular}
\caption{Left: Maya diagram  $(6,5,2,1,-2,-3)$ and its corresponding punctured young diagram. Right: The Maya diagram $(9,8,5,4,2)$ is the canonical representative of the shape class of previous diagram.  }
\label{fig:maya-partition}
\end{figure}

In accordance with our definition, we associate to every Maya diagram a
punctured Young diagram, which we define to be a Young diagram equiped with
an origin. This last should be regarded as marked point on the right border
of the Young diagram that separates the segments labelled by positive
integers from the segments labelled by strictly negative integers. An example can be seen in Figure \ref{fig:maya-partition}.

Supressing the origin, we recover a standard Young diagram, which we view as
an equivalence class of punctured Young diagrams and which we name a 
{\it shape class}. As a canonical representative in every shape class we take
the punctured Young diagram where the origin is located at the leftmost
point of the topmost horizontal half-line. The corresponding Maya diagram is
strictly positive. 

\subsection{Maya diagrams and partitions}

Every positive Maya diagram, i.e. every tuple of natural numbers $N_{m}=\left(
n_{1},...,n_{m}\right) $, is in one-to-one correspondence with an integer
partition \cite{Sergeev2014,GGM}%
\begin{equation}
\left\vert \lambda \right\vert =\lambda _{1}+...+\lambda _{m},
\end{equation}
where the $l(\lambda )=m$ parts of the partition $\lambda =\left( \lambda
_{1},...,\lambda _{m}\right) $ are given by

\begin{equation}
\lambda _{i}=n_{i}-m+i,
\end{equation}%
and%
\begin{equation}
\left\vert \lambda \right\vert =\sum_{i=1}^{m}n_{i}-m(m-1)/2.
\end{equation}

For the example in Figure \ref{fig:maya-partition}, we have  $N_{5}=\left( 9,8,5,4,2\right) $, the
corresponding partition is $\lambda =\left( 5^{2},3^{2},2\right) $ and $%
\left\vert \lambda \right\vert =18$ (as usual we let $\lambda _{i}^{k}$
denote $\lambda _{i}$ repeated $k$ times).

The $\lambda _{i}$ are the lengths of the columns of the corresponding Young
diagram starting from the lower left corner and $\lambda $ is a partition of
the integer $|\lambda|$.

Note that, contrary to the usual convention, we permit a partition to
terminate in a string of zeros. The reduced form $\widetilde{%
\lambda }$ of such a partition is obtained by suppressing this last string
of zeros. This corresponds to fixing the origin at the endpoint of the last
vertical segment on the boundary of the corresponding Young diagram. The
tuple corresponding to $\widetilde{\lambda }$ consists of strictly positive
integers; this last property characterizes the class of reduced partitions.
A reduced partition is therefore characteristic of the shape class and its
corresponding Young diagram.

To every reduced partition we can associate a conjugate partition $%
\overline{\lambda }$ \cite{macdonald,andrews} of the same integer $\left\vert
\lambda \right\vert $, whose elements $\overline{\lambda }_{j}$ are the
lengths of the rows (starting from the top) of the Young diagram associated
to $\lambda $. The conjugate partition can be visualized as the
partition corresponding to the transpose (with respect to the main diagonal
starting from the upperleft corner) of the original Young diagram. Formally,
we have

\begin{equation}
l(\overline{\lambda })=\lambda _{1}\text{,}\quad l(\lambda )=\overline{%
\lambda }_{1}  \label{conj1}
\end{equation}%
and 
\begin{equation}
\overline{\lambda }_{j}=\#\left\{ i:\lambda _{i}\geq j\right\} =\#\left\{
i:n_{i}\geq j-i+l(\lambda )\right\} .  \label{conj2}
\end{equation}

An alternative description of an irreductible partition, i.e. a shape class,
is via its Durfee symbol \cite{macdonald,andrews}. For that, we define the 
\textit{Durfee square} of the partition to be the largest square contained
in the Young diagram whose diagonal coincides with the main diagonal of the
diagram.

Letting $r\times r$ denote the size of the Durfee square, we define the
Durfee symbol to be the pair of partitions\bigskip\ 
\begin{equation}
\left[ \lambda ^{\prime }\mid \overline{\lambda }^{\prime }\right] =\left[
\lambda _{1}^{\prime },...,\lambda _{r}^{\prime }\mid \overline{\lambda }%
_{1}^{\prime },...,\overline{\lambda }_{r}^{\prime }\right] ,
\end{equation}%
where 
\begin{equation}
\lambda _{i}^{\prime }=\lambda _{i}-r,\quad \overline{\lambda }_{i}^{\prime
}=\overline{\lambda }_{i}-r.
\end{equation}%
Thus, the parts of $\lambda ^{\prime }$ are the length of the columns below
the Durfee square and the parts of $\overline{\lambda }^{\prime }$ are
length of the rows to the right of the Durfee square. In this way the
partition $\lambda ^{\prime }$ is a truncation of the original and $%
\overline{\lambda }^{\prime }$ is a truncation the conjugate partition.

\begin{figure}[ht]
\centering
\includegraphics{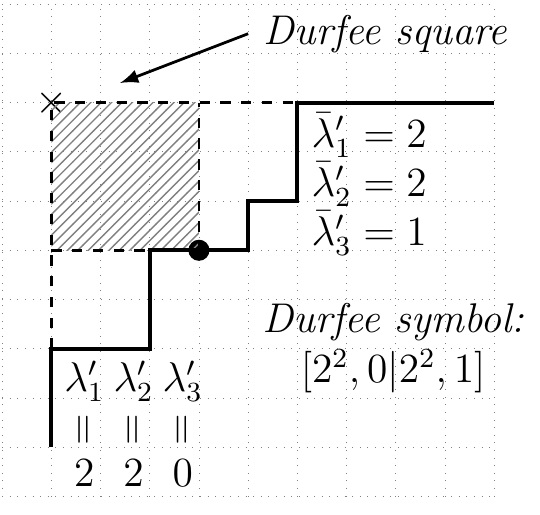}
\caption{Durfee square of the partition 
$\lambda =( 5^{2},3^{2},2)$ and corresponding Durfee symbol $\left[ 2^{2},0\mid 2^{2},1\right]$.}
\end{figure}

\subsection{Punctured Young diagrams and Durfee symbols}

The concept of Durfee square and Durfee symbol can be generalized to
describe an arbitrary punctured Young diagram, and correspondingly, every
Maya diagram. A Durfee rectangle is the $k\times \overline{k}$ rectangle
whose opposite vertices are the origin of a given punctured Young diagram
and the upper left corner of the diagram, and whose left side of length $k$
and upper side of length $\overline{k}$ lie on the left and upper border of
the corresponding Young diagram \cite{andrews1}. In this way we are able to
completely characterize a punctured Young diagram as the corresponding $%
k\times \overline{k}$ \ Durfee symbol, which we define as 
\begin{equation}
\left[ d\mid \overline{d}\right] _{k\times \overline{k}}=\left[
d_{1},...,d_{k}\mid \overline{d}_{1},...,\overline{d}_{\overline{k}}\right] ,
\label{ds}
\end{equation}%
where the 
\begin{equation}
d_{i}=\lambda _{i}-\overline{k}=n_{i}-k+i,  \label{corresp1}
\end{equation}%
stand for the heights of the truncated columns below the Durfee square and
where the 
\begin{equation}
\overline{d}_{j}=\overline{\lambda }_{j}-k=-\left( \overline{n}_{j}+1\right)
-\overline{k}+j,  \label{corresp2}
\end{equation}%
represent the lengths of the truncated rows at the right of the Durfee
rectangle. Note that

\begin{equation}
\begin{array}{l}
l(\lambda )=\left( k-\overline{k}\right) -\overline{n}_{1}=\overline{d}_{1}+k
\\ 
l(\overline{\lambda })=n_{1}+1-\left( k-\overline{k}\right) =d_{1}+\overline{%
k}.%
\end{array}
\label{l(lambda)}
\end{equation}

The Durfee rectangle can be considered as the rectangular Young diagram
corresponding to the partition

\begin{equation}
R=\left( \overline{k}^{\;k}\right) .  \label{R}
\end{equation}

From this point of view, $d$ and $\overline{d}$ are respectively the
partition and conjugate partition associated to each of the two connected
components of the diagram $\lambda - R$, (see \cite{andrews,macdonald}).

We define an irreducible Durfee symbol to be one for which $\overline{k}\leq
l(\overline{\lambda })$ and $k\leq l(\lambda )$, i.e. for which the punctured
Young diagram has an origin between its upper-right and lower left
corners. With the exception of the two extremal cases of $k=l(\overline{%
\lambda }),\overline{k}=0$ and $k=0,\overline{k}=l(\lambda )$, the
irreducible Durfee symbols are associated to non-flat rectangles. Thus, the
total number of such symbols is 
\begin{equation}
l(\lambda )+l(\overline{\lambda })+1=n_{1}-\overline{n}_{1}+2=d_{1}+%
\overline{d}_{1}+k+\overline{k}+1.
\end{equation}

\begin{figure}[ht]
\centering
\includegraphics{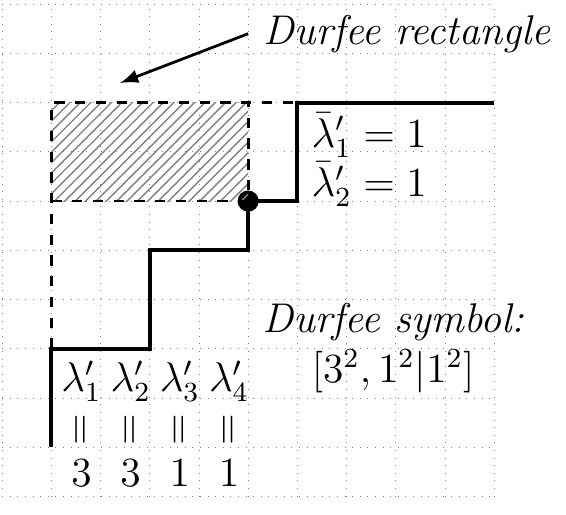}
\caption{Punctured Young diagram associated to the Maya diagram $(6,5,2,1,-2,-3)$, corresponding to a $2\times4$ Durfee
rectangle and Durfee symbol $[3^2,1^2|1^2]$.}
\end{figure}

If we fix the origin at the upper right corner of the Young diagram we have
the representative of the shape class. The associated Durfee rectangle is
horizontally flat and the corresponding generalized Durfee symbol is ($k=m$) 
$\left[ \lambda _{1},...,\lambda _{m}\mid \varnothing \right] $.

If we fix the origin at the lower left corner of the Young diagram, the
associated Durfee rectangle is vertically flat and the corresponding
generalized Durfee symbol reduces to $\left[ \varnothing \mid \overline{%
\lambda }_{1},...,\overline{\lambda }_{m}\right] $ $(\overline{k}=m).$

\vskip0.2cm
\subsection{Shifting rules}
\vskip0.2cm
Starting from the punctured Young diagram corresponding to $\left[ d\mid 
\overline{d}\right] $, a one-step displacement of the origin along the
border of the Young diagram changes the generalized Durfee symbol in a
simple way. The shifting rules are the following
\vskip0.2cm
\noindent\textbf{Left-descending shift}
\vskip0.2cm
We move of one step in the left (descending) direction. We denote by $_{1}\left[
d\mid \overline{d}\right] $ the resulting Durfee symbol

\begin{eqnarray}
d_{k}>0:&&\qquad _{1}\left[ d\mid \overline{d}\right] =\left[ d_{1}-1,...,d_{k}-1\mid 
\overline{d}_{1},...,\overline{d}_{\overline{k}},0\right] .  \label{shiftd}\\
d_{k}>0:&&\qquad_{1}\left[ d\mid \overline{d}\right] =\left[ d_{1},...,d_{k-1}\mid \overline{%
d}_{1}+1,...,\overline{d}_{\overline{k}}+1\right] .  \label{shiftd2}
\end{eqnarray}
From equations (\ref{corresp1}) and (\ref{corresp2}), we deduce the corresponding
updating rules for the associated Maya diagram $N_{m}$ :

\begin{eqnarray}
n_{k}>0:&&\qquad N_{m}\rightarrow \text{ }_{1}\left( N_{m}\right) =\left(
n_{1}-1,...,n_{k}-1,\left( -1\right) ,\overline{n}_{1}-1,...,\overline{n}_{%
\overline{k}}-1\right) .  \label{shiftn}\\
n_{k}>0:&&\qquad N_{m}\rightarrow \text{ }_{1}\left( N_{m}\right) =\left(
n_{1}-1,...,n_{k-1}-1,\overline{n}_{1}-1,...,\overline{n}_{\overline{k}%
}-1\right) .  \label{shiftn2}
\end{eqnarray}

\vskip0.2cm
\noindent\textbf{Right-ascending shift}
\vskip0.2cm

Starting from the punctured Young diagram corresponding to $\left[ d\mid 
\overline{d}\right] $, a one-step displacement of the origin along the
border of the Young diagram changes the generalized Durfee symbol in a
simple way. The unit left and right shifting rules are, respectively, the
following 
\begin{eqnarray}
_{1}\left[ d\mid \overline{d}\right] &=&%
\begin{cases}
\left[ d_{1}-1,...,d_{k}-1\mid \overline{d}_{1},...,\overline{d}_{\overline{k%
}},0\right] & \text{ if }d_{k}>0 \\ 
\left[ d_{1},...,d_{k-1}\mid \overline{d}_{1}+1,...,\overline{d}_{\overline{k%
}}+1\right] & \text{ if }d_{k}=0.%
\end{cases}
\label{shiftd3} \\
\left[ d\mid \overline{d}\right] _{1} &=&%
\begin{cases}
\left[ d_{1},...,d_{k},0\mid \overline{d}_{1}-1,...,\overline{d}_{\overline{k%
}}-1\right] & \text{ if }\overline{d}_{\overline{k}}>0 \\ 
\left[ d_{1}+1,...,d_{k}+1\mid \overline{d}_{1},...,\overline{d}_{\overline{k%
}-1}\right] & \text{ if }\overline{d}_{\overline{k}}=0.%
\end{cases}%
\end{eqnarray}

At the level of the Maya diagram, the above transformation corresponds to a
displacement of the origin one unit to the left and right, respectively. Let 
$N_{m}=\left( n_{1},...,n_{k},\overline{n}_{\overline{k}},...,\overline{n}%
_{1}\right) $ be the corresponding Maya diagram. From Eq \eqref{corresp1}
and Eq \eqref{corresp2}, we deduce the corresponding changes for the
associated Maya diagram : 
\begin{eqnarray}
_{1}\left( N_{m}\right) &=&%
\begin{cases}
\left( n_{1}-1,...,n_{k}-1,-1,\overline{n}_{1}-1,...,\overline{n}_{\overline{%
k}}-1\right) & \text{ if }n_{k}>0 \\ 
\left( n_{1}-1,...,n_{k-1}-1,\overline{n}_{1}-1,...,\overline{n}_{\overline{k%
}}-1\right) & \text{ if }n_{k}=0.%
\end{cases}
\label{shiftn3} \\
\left( N_{m}\right) _{1} &=&%
\begin{cases}
\left( n_{1}+1,...,n_{k}+1,0,\overline{n}_{1}+1,...,\overline{n}_{\overline{k%
}}+1\right) & \text{ if }\overline{n}_{\overline{k}}<-1 \\ 
\left( n_{1}+1,...,n_{k}+1,\overline{n}_{1}+1,...,\overline{n}_{\overline{k}%
-1}+1\right) & \text{ if }\overline{n}_{\overline{k}}=-1.%
\end{cases}%
\end{eqnarray}

In the considered example, $\left[ d\mid \overline{d}\right] =\left[
3^{2},1^{2}\mid 1^{2}\right] $, we have $11$ distinct irreducible Durfee
symbols in the associated shape class which are related by the following
shift connections (the right and left shifts are represented by right and
left arrows respectively)

\begin{eqnarray}
\left[ 5^{2},3^{2},2\mid \varnothing \right] &\rightleftarrows &\left[
4^{2},2^{2},1\mid 0\right] \rightleftarrows \left[ 3^{2},1^{2},0\mid 0^{2}%
\right] \rightleftarrows \left[ 3^{2},1^{2}\mid 1^{2}\right] \\
&\rightleftarrows &\left[ 2^{2},0^{2}\mid 1^{2},0\right] \rightleftarrows %
\left[ 2^{2},0\mid 2^{2},1\right] \rightleftarrows \left[ 2^{2}\mid 3^{2},2%
\right]  \notag \\
&\rightleftarrows &\left[ 1^{2}\mid 3^{2},2,0\right] \rightleftarrows \left[
0^{2}\mid 3^{2},2,0^{2}\right] \rightleftarrows \left[ 0\mid 4^{2},3,1^{2}%
\right] \rightleftarrows \left[ \varnothing \mid 5^{2},4,2^{2}\right] . 
\notag
\end{eqnarray}

For the corresponding Maya diagrams this becomes

\begin{eqnarray*}
\left( 9,8,5,4,2\right) &\rightleftarrows &\left( 8,7,4,3,1,-1\right)
\rightleftarrows \left( 7,6,3,2,0,-1,-2\right) \rightleftarrows \left(
6,5,2,1,-2,-3\right) \\
&\rightleftarrows &\left( 5,4,1,0,-1,-3,-4\right) \rightleftarrows \left(
4,3,0,-2,-4,-5\right) \rightleftarrows \left( 3,2,-3,-5,-6\right)  \notag \\
&\rightleftarrows &\left( 2,1,-1,-4,-6,-7\right) \rightleftarrows \left(
1,0,-1,-2,-5,-7,-8\right) \rightleftarrows \left( 0,-2,-3,-6,-8,-9\right)\\
& \rightleftarrows & \left( -3,-4,-7,-9,-10\right) .  \notag
\end{eqnarray*}

%%%%%%%%%%%%%%%%%%%%%%%%%%%%%%%%%%%%%%%%%%%%%%%%%%%%%%%%%%%%%%%%%%%
\section{Darboux Transformations}\label{sec:DT}
%%%%%%%%%%%%%%%%%%%%%%%%%%%%%%%%%%%%%%%%%%%%%%%%%%%%%%%%%%%%%%%%%%%

In this Section, and throughout the paper, the symbol $\propto$ will be heavily used  to denote proportionality up to a constant factor, i.e.
\[f(x)\propto g(x) \Leftrightarrow f(x)=C\cdot g(x), \quad C\in\mathbb R.\]

\subsection{One-step Darboux transformations}

We consider a one dimensional Hamiltonian $\widehat{H}=-d^{2}/dx^{2}+V(x),\
x\in I\subset R$ and the associated Schr\"{o}dinger equation%
\begin{equation}
\psi ^{\prime \prime }(x)+\left( E-V(x)\right) \psi (x)=0.  \label{EdS}
\end{equation}%
For every value of the eigenvalue $E$, there is a two dimensional space of
formal eigenfunctions $\psi (x)$ that satisfy the above equation. In the
following we also suppose that, with Dirichlet boundary conditions on $I$, $%
\widehat{H}$ admits a discrete spectrum of energies and eigenstates of the $%
\left( E_{n},\psi _{n}\right) _{n\in \left\{ 0,...,n_{\max }\right\} \mathbb{%
\ \subseteq N}}$ where, without loss of generality, we can always suppose
that the ground level of $\widehat{H}$ is zero: $E_{0}=0$. We consider $E$
as depending on an auxiliary spectral parameter $\lambda $, which identifies
to the quantum number $n$ for the energy levels. For a general value of $%
\lambda $, $\psi _{\lambda }$ means any solution of (\ref{EdS}) for the
eigenvalue $E_{\lambda }$ and for $\lambda =n$, it identifies the
eigenstate $\psi _{n}$.

Any formal eigenfunction $\psi _{\nu }(x)$ of $V(x)$ (i.e. $\widehat{H}$) can
be used as a seed function for a Darboux transformation (DT) $A(\psi _{\nu
}) $ which associates to the potential $V(x)$ a modified potential

\begin{equation}
V(x)\overset{A(\psi _{\nu })}{\longrightarrow }V^{\left( \nu \right)
}(x)=V(x)-2\left( \log \psi _{\nu }(x)\right) ^{\prime \prime },
\label{pottrans}
\end{equation}%
that we call an extension of $V(x)$. The formal eigenfunction of $V^{\left(
\nu \right) }$ associated to the spectral parameter $E_{\lambda }$ is given
by the Darboux-Crum formulas

{%
\begin{equation}
\psi _{\lambda }^{\left( \nu \right) }(x)\propto
\begin{cases}
\frac{W\left( \psi _{\nu },\psi _{\lambda }\mid x\right) }{\psi _{\nu }(x)}
& \text{ if }\lambda \neq \nu , \\ 
\frac{1}{\psi _{\nu }(x)}, & \text{ otherwise,}%
\end{cases}
\label{foDBTwronsk}
\end{equation}%
}where $W\left( y_{1},...,y_{m}\mid x\right) $ denotes the Wronskian of the
family of functions $y_{1},...,y_{m}$

\begin{equation}
W\left( y_{1},...,y_{m}\mid x\right) =\left\vert 
\begin{array}{ccc}
y_{1}\left( x\right) & ... & y_{m}\left( x\right) \\ 
... & \ddots & ... \\ 
y_{1}^{\left( m-1\right) }\left( x\right) & ... & y_{m}^{\left( m-1\right)
}\left( x\right)%
\end{array}%
\right\vert .  \label{wronskien}
\end{equation}

\subsection{Chains of Darboux transformations}

At the formal level, the  procedure can be straightforwardly iterated and a chain of 
$m$\ Darboux transformations can be completely characterized by the tuple $\left( \nu
_{1},...,\nu _{m}\right) $ of spectral indices of the successive seed
functions used in the chain. We note such a $m$-tuple by a capital letter $%
N_{m}$, where the index $m$ indicates the length of the chain $N_{m}=\left(
\nu _{1},...,\nu _{m}\right) $. $\psi _{\lambda }^{\left( N_{m}\right) }$ is
the formal eigenfunction associated to the eigenvalue $E_{\lambda }$ of the
potential $V^{\left( N_{m}\right) }(x)$. A chain is non-degenerate if all
the spectral indices $\nu _{i}$ of the chain $N_{m}$ are distinct and is
degenerate if some of them are repeated. In the rest of the paper, we use
the following notation

\begin{equation}
\left\{ 
\begin{array}{l}
\left( N_{m},\lambda \right) =\left( \nu _{1},...,\nu _{m},\lambda \right)
\\ 
\left( N_{m}+k\right) =\left( \nu _{1}+k,...,\nu _{m}+k\right).%
\end{array}%
\right.
\end{equation}

For non-degenerate chains, the extended potentials and their eigenfunctions
can be expressed via Wronskians of eigenfunctions of the initial potential 
\cite{crum,matveev}, known as \textit{Crum's formulas}. Assuming that all $\nu _{j}$ and $\lambda$ are distinct, we have

\begin{equation}
\psi _{\lambda }^{\left( N_{m}\right) }(x)=\frac{W^{\left( N_{m},\lambda
\right) }\left( x\right) }{W^{\left( N_{m}\right) }\left( x\right) },
\label{etats n3}
\end{equation}%
\begin{equation}
V^{\left( N_{m}\right) }(x)=V(x)-2\left( \log W^{\left( N_{m}\right) }\left(
x\right) \right) ^{\prime \prime },  \label{potnstep2}
\end{equation}%
where $W^{\left( N_{m}\right) }\left( x\right) =W\left( \psi _{\nu
_{1}},...,\psi _{\nu _{m}}\mid x\right) $.

We have shown \cite{GGM} that these formulas can be extended to degenerate
chains if we adopt the convention to suppress any pair of repeated index in
the set of spectral indices associated to the chain, namely

\begin{equation}
\psi _{\nu _{i}}^{\left( N_{m}\right) }(x)=\frac{W^{\left( \nu _{1},...,\nu
_{i-1},\nu _{i+1},...,\nu _{m}\right) }\left( x\right) }{W^{\left(
N_{m}\right) }\left( x\right) },  \label{etats n3deg}
\end{equation}%
which extends the second equality in (\ref{foDBTwronsk}) and

\begin{equation}
V^{\left( N_{m},\nu _{i}\right) }(x)=V(x)-2\left( \log W^{\left( \nu
_{1},...,\nu _{i-1},\nu _{i+1},...,\nu _{m}\right) }\left( x\right) \right)
^{\prime \prime }.  \label{potnstep2deg}
\end{equation}

%%%%%%%%%%%%%%%%%%%%%%%%%%%%%%%%%%%%%%%%%%%%%%%%%%%%%%%%%%%%%%%%%%%
\section{Translationally Shape Invariant Potentials}\label{sec:TSI}
%%%%%%%%%%%%%%%%%%%%%%%%%%%%%%%%%%%%%%%%%%%%%%%%%%%%%%%%%%%%%%%%%%%

Consider a potential $V(x;\alpha )$ which depends upon a (multi)parameter $%
\alpha =\left( \alpha ^{\left( 1\right) },...,\alpha ^{\left( M\right)
}\right) \in \mathbb{R}^{M}$ and with a (finite or infinite) bound state
spectrum $\left( E_{n},\psi _{n}\right) _{n\geq 0}$, the ground level being
supposed to be zero: $E_{0}(\alpha )=0$.

Any potential of the same functional form as $V(x;\alpha )$, i.e. of the form $%
V(x;\beta )+ K$, where $K$ is a constant,  will be said of the same \textit{shape class} as $V(x;\alpha )$. This potential will be said \textit{shape invariant} if its image by the Darboux transformation $A(\psi
_{0})$ with seed function equal to the ground state \cite{cooper,Dutt,gendenshtein}, is in the same shape class as the initial potential, i.e. if 
\begin{equation}
V^{\left( 0\right) }(x;\alpha )=V(x;f(\alpha ))+R(\alpha ),  \label{SI}
\end{equation}%
$R\left( \alpha \right) \in \mathbb{R}$ and $f(\alpha )\in \mathbb{R}^{M}$
being two given functions of $\alpha $. The corresponding Darboux transformation is usually
called a supersymmetric quantum mechanical (SUSY QM) partnership. The bound
state spectrum of such a shape invariant potential is then given by ($n\geq
1 $)

\begin{equation}
\left\{ 
\begin{array}{l}
E_{n}(\alpha )=\sum\limits_{k=0}^{n-1}R(\alpha
_{k})=\sum\limits_{k=0}^{n-1}E_{1}(\alpha _{k}) \\ 
\psi _{n}(x;\alpha )\propto \widehat{A}^{+}(\alpha )...\widehat{A}^{+}(\alpha
_{n-1})\psi _{0}(x;\alpha _{n}),%
\end{array}%
\right.  \label{spectreSIP}
\end{equation}%
where $\alpha _{k}=f^{\left( k\right) }(\alpha )=\overset{\text{k times}}{%
\overbrace{f\circ ...\circ f}}(\alpha )$ and $\widehat{A}^{+}(\alpha )=-%
\frac{d}{dx}-\left( \log \left( \psi _{0}(x;\alpha )\right) \right) ^{\prime
\prime }$.

When $f$ is a simple translation $f(\alpha )=\alpha +\varepsilon ,\
\varepsilon =\left( \varepsilon ^{\left( 1\right) },...,\varepsilon ^{\left(
M\right) }\right) $, $\varepsilon ^{\left( k\right) }=\pm 1$ or $0$, $V$ is
said to be \textit{translationally shape invariant} (TSI) and we call it a \textit{%
TSI Potential}. For all the known TSI potentials, we have $M\leq 2$.

For the TSI potentials, the \textit{dispersion relation} (energy as a function of
the quantum number) is  explicitly known and it has the general form
\begin{equation}
E_{n}(\alpha )=g(\alpha _{n})-g(\alpha ),  \label{g}
\end{equation}%
where  $g(\alpha)$ is a polynomial at most quadratic in $\alpha 
$. The dispersion relation can be extended analytically to non integer
values of the quantum number and we have in particular
\begin{equation}
E_{-\left( n+1\right) }(\alpha )=\sum\limits_{k=0}^{n}E_{-1}(\alpha
_{k}),\quad n\geq 0.  \label{Econj}
\end{equation}

A function is said to be \textit{quasi-polynomial} if it can be expressed,
up to a gauge factor, as an orthogonal polynomial in a suitable variable its
logarithmic derivative being rational in this variable. The eigenstates $\psi _{n}$ of a TSI potential  are always quasi-polynomial.

The set of TSI potentials contains in particular all the potentials classicaly known
to be exactly solvable, i.e. for which we know explicitly the dispersion
relation and whose the eigenfunctions can be expressed in closed analytical
form in terms of elementary transcendental functions, namely the harmonic,
isotonic, Morse, Kepler-Coulomb, Eckart, Darboux-P\"{o}schl-Teller
hyperbolic and trigonometric and Rosen-Morse hyperbolic and trigonometric
potentials. 

%These last are called \textbf{primary TSIP (PTSIP)}. Three of
%them are confining, the \textbf{harmonic oscillator (}HO), the \textbf{%
%isotonic oscillator (IO)} and the \textbf{trigonometric Darboux-P\"{o}%
%schl-Teller potential (tDPT)} potentials, defined respectively on the full
%real line ($I=%
%%TCIMACRO{\U{211d} }%
%%BeginExpansion
%\mathbb{R}
%%EndExpansion
%$), on the positive half line ($I=%
%%TCIMACRO{\U{211d} }%
%%BeginExpansion
%\mathbb{R}
%%EndExpansion
%^{+\ast }$) and on a finite interval ($I=\left] 0,\pi /2\right[ $) and which
%all have an infinite set of bound states $\left( E_{n},\psi _{n}\right)
%_{n\geq 0}$ constituting a complete family on $L^{2}(I)$.

In the rest of the paper, we concentrate on the isotonic oscillator and  trigonometric Darboux-P\"{o}schl-Teller (tDPT) potentials, whose 
eigenfunctions are expressible in terms of Laguerre and
Jacobi polynomials respectively. The dispersion relation $g(\alpha )$ is linear for the
isotonic oscillator and quadratic for the tDPT potential.

A common feature of these confining TSI potentials is the existence of a set of
symmetries $\Gamma _{i}$ which act on the parameters and which are
covariance transformations for the considered potential

\begin{equation}
\left\{ 
\begin{array}{l}
\alpha \overset{\Gamma _{i}}{\rightarrow }\Gamma _{i}(\alpha ) \\ 
V(x;\alpha )\overset{\Gamma _{i}}{\rightarrow }V(x;\Gamma _{i}(\alpha
))=V(x;\alpha )+\delta _{i}\left( \alpha \right) .%
\end{array}%
\right.  \label{sym}
\end{equation}

In the space of parameters, $\Gamma _{i}$ correspond to the reflections
with respect to the coordinate axis. For the isotonic oscillator and tDPT potentials, the parameter space is two dimensional and we have three distinct
symmetries $\Gamma _{1}$, $\Gamma _{2}$ and $\Gamma _{3}=\Gamma _{1}\circ $ $%
\Gamma _{2}$.

A formal eigenfunction diverging at both endpoints of the definition
interval will be said of $(\infty,\infty)$  type. A
formal eigenfuntion diverging at one or the other extremity of the
definition interval will be said of $( 0,\infty)$ or
$(\infty ,0)$ type and a formal eigenfuntion tending
to zero at both endpoints of the definition interval will be said of $(0,0)$ type, which includes the $L^2$ eigenfunctions.

\subsection{Pseudo spectra of the trigonometric Darboux-P\"{o}schl-Teller potential}

The trigonometric Darboux-P\"{o}schl-Teller (tDPT) potential  is defined on the interval $(0,\pi /2)$ by

\begin{equation}
V(x;\alpha ,\beta )=\frac{(\alpha +1/2)(\alpha -1/2)}{\sin ^{2}x}+\frac{%
(\beta +1/2)(\beta -1/2)}{\cos ^{2}x}-(\alpha +\beta +1)^{2},\ \left\vert
\alpha \right\vert ,\left\vert \beta \right\vert >1/2.  \label{tDPT}
\end{equation}

With Dirichlet boundary conditions at $0$ and $\pi /2$ and in the case $%
\alpha ,\beta >1/2$, it has the following spectrum 
\begin{equation}
\left\{ 
\begin{array}{l}
\text{ }E_{n}\left( \alpha ,\beta \right) =(\alpha _{n}+\beta
_{n}+1)^{2}-(\alpha +\beta +1)^{2}=4n(\alpha +\beta +1+n) \\ 
\\ 
\psi _{n}\left( x;\alpha ,\beta \right) =\left( 1-z\right) ^{\left( \alpha
+1/2\right) /2}\left( 1+z\right) ^{\left( \beta +1/2\right) /2}\mathit{P}%
_{n}^{\left( \alpha ,\beta \right) }\left( z\right)%
\end{array}%
\right. ,\ n\in \mathbb{N},  \label{spec tDPT}
\end{equation}%
where $z=\cos 2x$ ($\psi _{n}$ is of $\left( 0,0\right) $ type on $\left[
0,\pi /2\right] $). The $\mathit{P}_{n}^{\left( \alpha ,\beta \right) }(z)$
are the usual Jacobi polynomials \cite{szego,magnus} and $(\alpha _{n},\beta
_{n})=(\alpha +n,\beta +n)$.

Note that $V(x;\alpha,\beta)$ is invariant under the parity reflection $x\rightarrow \pi /2-x$ ( which corresponds to $z\rightarrow -z$) combined with the permutation of the parameters ($\alpha
\leftrightarrow \beta $). This implies in a direct way the standard
reflection property of the Jacobi polynomials \cite{szego,magnus}:

\begin{equation}
\mathit{P}_{n}^{\left( \alpha ,\beta \right) }\left( -z\right) =\left(
-1\right) ^{n}\mathit{P}_{n}^{\left( \beta ,\alpha \right) }\left( z\right) .
\label{reflJac}
\end{equation}

The tDPT\ potential is a TSI potential, with $\left( \alpha ,\beta \right) \in 
%TCIMACRO{\U{211d} }%
%BeginExpansion
\mathbb{R}
%EndExpansion
^{2}$ and $\varepsilon =\left( +1,+1\right) $

\begin{equation}
\left\{ 
\begin{array}{l}
V^{\left( 0\right) }\left( x;\alpha ,\beta \right) =V\left( x;\alpha
_{1},\beta _{1}\right) +4(\alpha +\beta +2) \\ 
\psi _{n}^{\left( 0\right) }\left( x;\alpha ,\beta \right) \propto \psi
_{n-1}\left( x;\alpha _{1},\beta _{1}\right)%
\end{array}%
\right. .  \label{SI tDPT}
\end{equation}

It possesses three discrete parametric symmetries:

\begin{enumerate}
\item The $\Gamma _{1}$  symmetry which acts as

\begin{equation}
\left\{ 
\begin{array}{l}
\left( \alpha ,\beta \right) \overset{\Gamma _{1}}{\rightarrow }\left(
\alpha ,-\beta \right) \\ 
V(x;\alpha ,\beta )\overset{\Gamma _{1}}{\rightarrow }V(x;\alpha ,-\beta
)=V(x;\alpha ,\beta )+4\beta (\alpha +1),%
\end{array}%
\right.
\end{equation}%
and generates the \textit{conjugate shadow spectrum} of $V\left( x;\alpha ,\beta
\right) $ :
\begin{equation}
\left\{ 
\begin{array}{l}
E_{n}\left( \alpha ,\beta \right) \\ 
\psi _{n}(x;\alpha ,\beta )%
\end{array}%
\right. \overset{\Gamma _{1}}{\rightarrow }\left\{ 
\begin{array}{l}
E_{n-\beta }\left( \alpha ,\beta \right) =4\left( n-\beta \right) (n+\alpha
+1)=E_{-\left( n+1\right) -\alpha }\left( \alpha ,\beta \right) \\ 
\phi _{-\left( n+1\right) }\left( x;\alpha ,\beta \right) =\psi _{n}\left(
x;\alpha ,-\beta \right) =\left( 1-z\right) ^{\left( \alpha +1/2\right)
/2}\left( 1+z\right) ^{\left( -\beta +1/2\right) /2}\mathit{P}_{n}^{\left(
-\alpha ,\beta \right) }\left( z\right) ,%
\end{array}%
\right.  \label{conjshadtDPT}
\end{equation}
where 
\[ \phi _{-\left( n+1\right) }\left( x;\alpha ,\beta \right) =\psi
_{n-\beta }\left( x;\alpha ,\beta \right) =\psi _{-\left( n+1\right) -\alpha
}\left( x;\alpha ,\beta \right) \]
 is of $\left( 0,\infty \right) $ type on $\left[ 0,\pi /2\right] $. Note also that

\begin{eqnarray}
\Gamma _{1}\left( V(x;\alpha _{k},\beta _{l})\right) &=&V(x;\alpha _{k},\beta
_{-l})+4\beta _{-l}\left( \alpha _{k}+1\right) \\ 
\Gamma _{1}\left( \psi _{n}(x;\alpha _{k},\beta _{l})\right) &=&\phi
_{n}(x;\alpha _{k},\beta _{-l}),\quad n\in \mathbb{Z}
\end{eqnarray}

\item The $\Gamma _{2}$ symmetry which acts as

\begin{equation}
\left( \alpha ,\beta \right) \overset{\Gamma _{2}}{\rightarrow }\left(
-\alpha ,\beta \right) ,\left\{ 
\begin{array}{l}
\left( \alpha ,\beta \right) \overset{\Gamma _{2}}{\rightarrow }\left(
-\alpha ,\beta \right) \\ 
V(x;\alpha ,\beta )\overset{\Gamma _{2}}{\rightarrow }V(x;-\alpha ,\beta
)=V(x;\alpha ,\beta )+4\alpha (\beta +1)%
\end{array}%
\right.
\end{equation}%
and generates the \textit{shadow spectrum} of $V\left( x;\alpha ,\beta \right) $:
\begin{equation}
\left\{ 
\begin{array}{l}
E_{n}\left( \alpha ,\beta \right) \\ 
\psi _{n}(x;\alpha ,\beta )%
\end{array}%
\right. \overset{\Gamma _{2}}{\rightarrow }\left\{ 
\begin{array}{l}
E_{n-\alpha }\left( \alpha ,\beta \right) =4\left( n-\alpha \right) (n+\beta
+1) \\ 
\phi _{n}\left( x;\alpha ,\beta \right) =\psi _{n}\left( x;-\alpha ,\beta
\right) =\left( 1-z\right) ^{\left( -\alpha +1/2\right) /2}\left( 1+z\right)
^{\left( \beta +1/2\right) /2}\mathit{P}_{n}^{\left( -\alpha ,\beta \right)
}\left( z\right) ,%
\end{array}%
\right.  \label{shadtDPT}
\end{equation}%
where $\phi _{n}\left( x;\alpha ,\beta \right) =\psi _{n-\alpha }\left(
x;\alpha ,\beta \right) $ is of $\left( \infty ,0\right) $ type on $\left[
0,\pi /2\right] $. We have also

\begin{equation}
\left\{ 
\begin{array}{l}
\Gamma _{2}\left( V(x;\alpha _{k},\beta _{l})\right) =V(x;\alpha _{-k},\beta
_{l})+4\alpha _{-k}\left( \beta _{l}+1\right) \\ 
\Gamma _{2}\left( \psi _{n}(x;\alpha _{k},\beta _{l})\right) =\phi _{-\left(
n+1\right) }(x;\alpha _{-k},\beta _{l}),\quad n\in 
%TCIMACRO{\U{2124} }%
%BeginExpansion
\mathbb{Z}
%EndExpansion
.%
\end{array}%
\right.
\end{equation}

\item The $\Gamma _{3}=\Gamma _{1}\circ \Gamma _{2}$ symmetry which
acts as
\begin{equation}
\left( \alpha ,\beta \right) \overset{\Gamma _{3}}{\rightarrow }\left(
-\alpha ,-\beta \right) ,\left\{ 
\begin{array}{l}
V(x;\alpha ,\beta )\overset{\Gamma _{3}}{\rightarrow }V(x;-\alpha ,-\beta
)=V(x;\alpha ,\beta )+4(\alpha +\beta ) \\ 
\psi _{n}(x;\alpha ,\beta )\overset{\Gamma _{3}}{\rightarrow }\psi _{-\left(
n+1\right) }\left( x;\alpha ,\beta \right) ,%
\end{array}%
\right.  \label{gamma3DPT}
\end{equation}%
and generates the \textit{conjugate spectrum} of $V\left( x;\alpha ,\beta
\right) $.
\begin{equation}
\left\{ 
\begin{array}{l}
E_{n}\left( \alpha ,\beta \right) \\ 
\psi _{n}(x;\alpha ,\beta )%
\end{array}%
\right. \overset{\Gamma _{3}}{\rightarrow }\left\{ 
\begin{array}{l}
E_{-\left( n+1\right) }\left( \alpha ,\beta \right) =-4\left( n+1\right)
(\alpha +\beta -n) \\ 
\psi _{-\left( n+1\right) }\left( x;\alpha ,\beta \right) =\psi
_{n}(x;-\alpha ,-\beta )=\left( 1-z\right) ^{\left( -\alpha +1/2\right)
/2}\left( 1+z\right) ^{\left( -\beta +1/2\right) /2}\mathit{P}_{n}^{\left(
-\alpha ,-\beta \right) }\left( z\right) ,%
\end{array}%
\right.  \label{conjtDPT}
\end{equation}%
where $\psi _{-\left( n+1\right) }\left( x;\alpha ,\beta \right) $ is of $%
\left( \infty ,\infty \right) $ type. Note also that

\begin{equation}
\left\{ 
\begin{array}{l}
\Gamma _{3}\left( V(x;\alpha _{k},\beta _{l})\right) =V(x;\alpha _{-k},\beta
_{-l})+4\left( \alpha _{-k}+\beta _{-l}\right) \\ 
\Gamma _{3}\left( \psi _{n}(x;\alpha _{k},\beta _{l})\right) =\psi _{-\left(
n+1\right) }(x;\alpha _{-k},\beta _{-l}),\quad n\in 
%TCIMACRO{\U{2124} }%
%BeginExpansion
\mathbb{Z}
%EndExpansion
.%
\end{array}%
\right.
\end{equation}

\end{enumerate}

The conjugate spectrum is obtained by a mirror symmetry (center at $\nu
=-1/2 $) on the values of the spectral parameter: $\nu \rightarrow -(\nu +1)$.
The tDPT potential presents in addition to the true spectrum, the three 
\textit{pseudo-spectra} defined above: conjugate ($E_{-\left( n+1\right) }$%
), shadow ($E_{n-\alpha }$) and conjugate shadow ($E_{n-\beta }$) spectra.
The union of the quasi-spectrum and of its conjugate will be called the 
\textit{extended spectrum}. The union of the shadow spectrum and of its
conjugate will be called the \textit{extended shadow spectrum}.

The extended shadow spectrum is thus obtained by shifting the integer values
of the spectral parameter associated to the spectrum, by $\beta $ and $%
\alpha $. Due to the structure of the dispersion relation, which induces the
identity $E_{n-\beta }=E_{-\left( n+1\right) -\alpha }$, it can also be
 obtained by shifting the values of the extended and conjugate
spectra by $\alpha $. Under this aspect, the conjugate shadow spectrum $%
E_{-\left( n+1\right) -\alpha }$ appears to be the shadow of the conjugate
spectrum $E_{-\left( n+1\right) }$. Starting from the usual eigenstates, the 
$\Gamma _{i}$ symmetries generate all the quasi-polynomial eigenfunctions of
the tDPT potential which are all contained in the extended or extended
shadow spectra.

If $\alpha \neq \beta $, and $\alpha $ and $\beta $ are not integers, the four
pseudo-spectra are distinct and in the rest of the paper we limit our study to this
non-degenerate case. For integer values of $\alpha $ or $\beta $ the pseudo-spectra can
coalesce, giving rise to interesting but specific phenomena \cite{GCB}.

\subsection{Pseudo spectra of the Isotonic oscillator}

The isotonic oscillator potential (with zero ground level $E_{0}=0$) is defined on the
positive half line $(0,+\infty) $ by
\begin{equation}
V\left( x;\omega ,\alpha \right) =\frac{\omega ^{2}}{4}x^{2}+\frac{\left(
\alpha +1/2\right) (\alpha -1/2)}{x^{2}}-\omega \left( \alpha +1\right)
,\quad \left\vert \alpha \right\vert >1/2.  \label{OI}
\end{equation}

If we impose Dirichlet boundary conditions at $0$ and infinity and assume $\alpha >1/2$ and $\alpha \notin \mathbb N$, it has the following spectrum ($z=\omega x^{2}/2$)%
\begin{equation}
\left\{ 
\begin{array}{l}
E_{n}\left( \omega \right) =2n\omega \\ 
\psi _{n}\left( x;\omega ,\alpha \right) =z^{\left( \alpha +1/2\right)
/2}e^{-z/2}\mathit{L}_{n}^{\alpha }\left( z\right)%
\end{array}%
\right. ,\quad n\geq 0,  \label{spec OI}
\end{equation}%
where $\psi _{n}$ is of $\left( 0,0\right) $ type on $[ 0,+\infty )$. The isotonic oscillator is a TSI potential, with $\left( \omega ,\alpha \right) \in 
%TCIMACRO{\U{211d} }%
%BeginExpansion
\mathbb{R}
%EndExpansion
^{2}$ and $\varepsilon =\left( 0,+1\right) $, since
\begin{equation}
\left\{ 
\begin{array}{l}
V^{\left( 0\right) }\left( x;\omega ,\alpha \right) =V\left( x;\omega
,\alpha _{1}\right) +2\omega \\ 
\psi _{n}^{\left( 0\right) }\left( x;\omega ,\alpha \right) \propto \psi
_{n-1}\left( x;\omega ,\alpha _{1}\right)%
\end{array}%
\right. ,  \label{SI IO}
\end{equation}%
with $\alpha _{n}=\alpha +n$. This potential also possesses three discrete parametric symmetries, which we describe below.

\begin{enumerate}
\item The $\Gamma _{1}$ symmetry which acts as
\begin{equation}
\left\{ 
\begin{array}{l}
\left( \omega ,\alpha \right) \overset{\Gamma _{1}}{\rightarrow }\left(
-\omega ,\alpha \right) \\ 
V(x;\omega ,\alpha )\overset{\Gamma _{1}}{\rightarrow }V(x;-\omega ,\alpha
)=V(x;\omega ,\alpha )+2\omega \left( \alpha +1\right) ,%
\end{array}%
\right.
\end{equation}%
and generates the \textit{conjugate shadow spectrum} of $V\left( x;\omega
,\alpha \right) $ :
\begin{equation}
\left\{ 
\begin{array}{l}
E_{n}\left( \omega \right) \\ 
\psi _{n}(x;\omega ,\alpha )%
\end{array}%
\right. \overset{\Gamma _{1}}{\rightarrow }\left\{ 
\begin{array}{l}
E_{-\left( n+1\right) -\alpha }\left( \omega \right) =-2\left( n+1+\alpha
\right) \omega <0 \\ 
\phi _{-(n+1)}(x;\omega ,\alpha )=\psi _{n}(x;-\omega ,\alpha )=z^{\left(
\alpha +1/2\right) /2}e^{z/2}\mathit{L}_{n}^{\alpha }\left( -z\right)%
\end{array}%
\right. ,\quad n\geq 0,  \label{conjshadOI}
\end{equation}%
where $\phi _{-(n+1)}(x;\omega ,\alpha )=\psi _{-(n+1)-\alpha }(x;\omega
,\alpha )$ is of $\left( 0,\infty \right) $ type on $[ 0,+\infty )
$. Note also that

\begin{equation}
\left\{ 
\begin{array}{l}
\Gamma _{1}\left( V(x;\omega ,\alpha _{k})\right) =V(x;\omega ,\alpha
_{k})+2\omega \left( \alpha _{k}+1\right) \\ 
\Gamma _{1}\left( \psi _{n}(x;\omega ,\alpha _{k})\right) =\phi _{-\left(
n+1\right) }(x;\omega ,\alpha _{k}),\quad n\in 
%TCIMACRO{\U{2124} }%
%BeginExpansion
\mathbb{Z}
%EndExpansion
.%
\end{array}%
\right.
\end{equation}

\item The $\Gamma _{2}$ symmetry which acts as

\begin{equation}
\left\{ 
\begin{array}{l}
\left( \omega ,\alpha \right) \overset{\Gamma _{2}}{\rightarrow }\left(
\omega ,-\alpha \right) \\ 
V(x;\omega ,\alpha )\overset{\Gamma _{2}}{\rightarrow }V(x;\omega ,-\alpha
)=V(x;\omega ,\alpha )+2\omega \alpha ,%
\end{array}%
\right.
\end{equation}%
and generates the \textit{shadow spectrum} of $V\left( x;\omega ,\alpha \right)
: $
\begin{equation}
\left\{ 
\begin{array}{l}
E_{n}\left( \omega \right) \\ 
\psi _{n}(x;\omega ,\alpha )%
\end{array}%
\right. \overset{\Gamma _{2}}{\rightarrow }\left\{ 
\begin{array}{l}
E_{n-\alpha }\left( \omega \right) =2\left( n-\alpha \right) \omega \\ 
\phi _{n}(x;\omega ,\alpha )=\psi _{n}(x;\omega ,-\alpha )=z^{\left( -\alpha
+1/2\right) /2}e^{-z/2}\mathit{L}_{n}^{-\alpha }\left( z\right)%
\end{array}%
\right. ,\quad n\geq 0,  \label{shadOI}
\end{equation}%
where $\phi _{n}(x;\omega ,\alpha )$ is of $\left( \infty ,0\right) $ type
on $[ 0,+\infty)$. We have in particular%
\begin{equation}
\phi _{-1}(x;\omega ,\alpha )=z^{\left( \alpha +1/2\right) /2}e^{z/2}=1/\phi
_{0}(x;\omega ,\alpha _{1}).
\end{equation}%
and more generally
\begin{equation}
\left\{ 
\begin{array}{l}
\Gamma _{2}\left( V(x;\omega ,\alpha _{k})\right) =V(x;\omega ,\alpha
_{-k})+2\omega \alpha _{-k} \\ 
\Gamma _{2}\left( \psi _{n}(x;\omega ,\alpha _{k})\right) =\phi
_{n}(x;\omega ,\alpha _{-k}),\quad n\in 
%TCIMACRO{\U{2124} }%
%BeginExpansion
\mathbb{Z}
%EndExpansion
.%
\end{array}%
\right.
\end{equation}

\item The $\Gamma _{3}=\Gamma _{1}\circ \Gamma _{2}$ symmetry which
acts as
\begin{equation}
\left\{ 
\begin{array}{l}
\left( \omega ,\alpha \right) \overset{\Gamma _{3}}{\rightarrow }\left(
-\omega ,-\alpha \right) \\ 
V(x;\omega ,\alpha )\overset{\Gamma _{3}}{\rightarrow }V(x;-\omega ,-\alpha
)=V(x;\omega ,\alpha )+2\omega ,%
\end{array}%
\right.
\end{equation}%
and generates the \textit{conjugate spectrum} of $V\left( x;\omega ,\alpha
\right) $ :
\begin{equation}
\left\{ 
\begin{array}{l}
E_{n}\left( \omega \right) \\ 
\psi _{n}(x;\omega ,\alpha )%
\end{array}%
\right. \overset{\Gamma _{3}}{\rightarrow }\left\{ 
\begin{array}{l}
E_{-\left( n+1\right) }\left( \omega \right) =-2\left( n+1\right) \omega <0
\\ 
\psi _{-\left( n+1\right) }(x;\omega ,\alpha )=\psi _{n}(x;-\omega ,-\alpha
)=z^{\left( -\alpha +1/2\right) /2}e^{z/2}\mathit{L}_{n}^{-\alpha }\left(
-z\right)%
\end{array}%
\right. ,\quad n\geq 0,  \label{conjOi}
\end{equation}%
where $\psi _{-\left( n+1\right) }(x;\omega ,\alpha )$ is of $\left( \infty
,\infty \right) $ type on $[ 0,+\infty ) $. In particular, we have
\begin{equation}
\psi _{-1}(x;\omega ,\alpha )=z^{\left( -\alpha +1/2\right)
/2}e^{z/2}=1/\psi _{0}(x;\omega ,\alpha _{-1}).  \label{psi-1}
\end{equation}
and more generally
\begin{equation}
\left\{ 
\begin{array}{l}
\Gamma _{3}\left( V(x;\omega ,\alpha _{k})\right) =V(x;\omega ,\alpha
_{-k})+2\omega \\ 
\Gamma _{3}\left( \psi _{n}(x;\omega ,\alpha _{k})\right) =\psi _{-\left(
n+1\right) }(x;\omega ,\alpha _{-k}),\quad n\in 
%TCIMACRO{\U{2124} }%
%BeginExpansion
\mathbb{Z}
%EndExpansion
.%
\end{array}%
\right.  \label{conjOi2}
\end{equation}

\end{enumerate}

Here again, the presence of the infinite centrifugal barrier at the origin
induces the existence of an extended shadow spectra, $\left( E_{n-\alpha
}\right) _{n\in 
%TCIMACRO{\U{2124} }%
%BeginExpansion
\mathbb{Z}
%EndExpansion
}$. It is obtained by shifting the integer values of the spectral parameter
associated to the extended spectrum $\left( E_{n}\right) _{n\in 
%TCIMACRO{\U{2124} }%
%BeginExpansion
\mathbb{Z}
%EndExpansion
}$ by $\alpha $. Note that, the range $-1/2<\alpha <1/2$ that we have
excluded previously is nevertheless interesting because then the potential
is weakly attractive at $z=0$. This attraction is indeed sufficiently weak
that we still get a self-adjoint eigenvalue problem but for which one has to
use non-Dirichlet boundary conditions.

\subsection{Mixed Darboux Chains for TSI potentials}

To produce the rational extensions of the considered potentials we have at
our disposal all the quasipolynomial eigenfunctions, that can be used as seed functions in Darboux transformations:
\begin{itemize}
\item the eigenstates and their conjugates which can be gathered into one family 
$\psi _{n}$, labelled by positive and negative integer values of the quantum
number $n\in \mathbb{Z}$ and which form the \textit{extended spectrum}.

\item the shadow and conjugate shadow eigenfunctions which can be gathered into
one family $\phi _{n}$, labelled by positive and negative integer values of
the quantum number $n\in \mathbb{Z}$ and which form the \textit{extended
shadow spectrum}.
\end{itemize}

Using $\psi _{n}$ and $\phi _{n},\ n\in \mathbb Z$ as seed functions for chains of Darboux transformations, the obtained potentials are rational
functions in the adapted variable. In the general case, the tuple of
spectral indices associated to the chain contains indices associated to the
extended spectrum $\left( n_{1},...,n_{m}\right) $ and indices associated to
the extended shadow spectrum $\left( l_{1},...,l_{r}\right) $ and the
corresponding chain of Darboux transformations will be called a \textit{mixed chain}. The chain is then characterized by a bi-tuple of spectral indices
\begin{equation}
N_{m}\otimes L_{r}=\left( n_{1},...,n_{m}\right) \otimes \left(
l_{1},...,l_{r}\right) .  \label{bi-tuple}
\end{equation}
If either $m=0$ or $r=0$, the bi-tuple reduces to a standard tuple $L_{r}$ or $%
N_{m} $ and the chain is called a \textit{single type chain}. A \textit{reducible chain} is a single type chain whose associated tuple
contains consecutive integers starting from $0$ in increasing order or
starting from $ -1 $ in decreasing order.

The bi-tuple $N_{m}\otimes L_{r}$ can be associated to a pair of Maya
diagrams, one for $N_{m}$ and the other one for $L_{r}$.

In the rest of the paper we will use the following notation: 
\begin{eqnarray*}
N_{m}&=&\left( n_{1},...,n_{k},\overline{n}_{\overline{k}},...,\overline{n}_{1}\right) \text{ with } k+\overline{k}=m.\quad  n_{i}>n_{i+1}\geq 0, \quad 0>\overline{n}_{i+1}>\overline{n}_{i},\\
L_{r}&=&\left( l_{1},...,l_{\kappa },\overline{l}_{\kappa },...,\overline{l}
_{1}\right),\text{ with }\kappa +\overline{\kappa }=r,\quad l_{i}>l_{i+1}\geq 0,\quad 0>%
\overline{l}_{i+1}>\overline{l}_{i}.
\end{eqnarray*}
The bi-tuple $N_{m}\otimes L_{r}$ can be equivalently associated to a pair
of punctured Young diagrams or to a direct product Durfee symbols (see equations
(\ref{corresp1}) and (\ref{corresp2}))
\begin{equation}
\left[ d\mid \overline{d}\right] \otimes \left[ \delta \mid \overline{\delta 
}\right] =\left[ d_{1},...,d_{k}\mid \overline{d}_{1},...,\overline{d}_{%
\overline{k}}\right] \otimes \left[ \delta _{1},...,\delta _{\kappa }\mid 
\overline{\delta }_{1},...,\overline{\delta }_{\overline{\kappa }}\right] ,
\label{PUC}
\end{equation}%
where $d_{i}=n_{i}-k+i$, $\delta _{i}=l_{i}-\kappa +i$\ and $\overline{d}%
_{j}=-\left( \overline{n}_{i}+1\right) -\overline{k}+i$, $\overline{\delta }%
_{j}=-\left( \overline{l}_{i}+1\right) -\overline{\kappa }+i$.

Refering to Koike \cite{koike} and Tsuda \cite{tsuda}, we call $\left[ d\mid 
\overline{d}\right] \otimes \left[ \delta \mid \overline{\delta }\right] $ a 
\textit{punctured universal character}, which can be thought of as the direct product of two punctured Young diagrams.

In \cite{GGM1}, using shape invariance arguments, in the case of mixed
chains containing only seed functions of the extended spectrum, we have
proven that every mixed chain is equivalent to a single type chain $%
N_{m}=\left( n_{1},...,n_{m}\right) \in 
%TCIMACRO{\U{2115} }%
%BeginExpansion
\mathbb{N}
%EndExpansion
^{m}$ of seed functions of the spectrum. In other words these type
of extensions can be gathered into shape classes of potentials whose
representative element is an extension with an associated Durfee symbol of
the form $\left[ \lambda _{1},...,\lambda _{m}\mid \varnothing \right] $.
Such a shape class is in one to one correspondence with the Young diagram of
the partition $\lambda $ which represents a shape class of punctured Young
diagrams. Using different\ approaches, Odake \cite{odake} has obtained
equivalent results for chains containing only seed functions of the extended
shadow spectrum and Takemura \cite{takemura} has extended this results for
general mixed chains for the isotonic oscillator and tDPT potentials.

As we will see, all these results are direct consequences of the
shape invariance and parametric symmetries of the primary potential. The
resulting equivalence formulas are just obtained by applying
shape invariance for its rationally extended descendants. In fact, all the rational extensions of the three primary confining TSI potentials and all the associated exceptional orthogonal polynomials, can be put in one to one correspondence with a punctured universal character.

The equivalence must be understood in the following manner at the level of potentials and Durfee symbols: shifting the origin in a punctured Young diagram modifies the chain of Darboux transformations but produces a potential in the same shape class. Equivalently, the pseudowronskians with orthogonal polynomial entries (which are, up to gauge factors, the Wronskians of the seed functions associated to the chains) corresponding to equivalent symbols will be identical up to an overall factor.

%%%%%%%%%%%%%%%%%%%%%%%%%%%%%%%%%%%%%%%%%%%%%%%%%%%%%%%%%%%%%%%%%%%
\section{Shape invariance and formulas for the derivatives of Laguerre and Jacobi polynomials}\label{sec:shapeinv}
%%%%%%%%%%%%%%%%%%%%%%%%%%%%%%%%%%%%%%%%%%%%%%%%%%%%%%%%%%%%%%%%%%%

\subsection{Isotonic oscillator and Laguerre polynomials}

In terms of punctured universal characters, the shape invariance property (\ref{SI IO}) for the isotonic oscillator becomes
\begin{equation}
\left\{ 
\begin{array}{l}
V^{\left[ 0\mid \varnothing \right] \otimes \left[ \varnothing \mid
\varnothing \right] }(x;\omega ,\alpha )=V^{\left( 0\right) \otimes \left(
\varnothing \right) }(x;\omega ,\alpha )=V(x;\omega ,\alpha
_{1})+E_{1}(\omega ) \\ 
\psi _{n}^{\left[ 0\mid \varnothing \right] \otimes \left[ \varnothing \mid
\varnothing \right] }(x;\omega ,\alpha )=\psi _{n}^{\left( 0\right) \otimes
\left( \varnothing \right) }(x;\omega ,\alpha )\propto \psi _{n-1}(x;\omega
,\alpha _{1})%
\end{array}%
\right. ,n\in 
%TCIMACRO{\U{2124} }%
%BeginExpansion
\mathbb{Z}
%EndExpansion
.  \label{SIIO1}
\end{equation}

In the particular case  $n=0$, we have by (\ref{conjOi}):
\[\psi _{0}^{\left[ 0\mid \varnothing %
\right] \otimes \left[ \varnothing \mid \varnothing \right] }(x;\omega
,\alpha )=1/\psi _{0}(x;\omega ,\alpha )\propto \psi _{-1}(x;\omega ,\alpha
_{1}).\]

Observe that $\phi _{n}^{\left[ 0\mid \varnothing \right] \otimes \left[
\varnothing \mid \varnothing \right] }(x;\omega ,\alpha )=\phi _{n}^{\left(
0\right) \otimes \left( \varnothing \right) }(x;\omega ,\alpha )$ is a
formal eigenfunction of $V^{\left[ 0\mid \varnothing \right] \otimes \left[
\varnothing \mid \varnothing \right] }(x;\omega ,\alpha )$ for the formal
eigenvalue $E_{n-\alpha }(\omega )$, that is, a formal eigenfunction of $%
V(x;\omega ,\alpha _{1})$ for the formal eigenvalue $E_{n-\alpha }(\omega
)-E_{1}(\omega )=E_{n-\alpha _{1}}(\omega )$. Since it is also
quasi-polynomial we necessarily have

\begin{equation}
\phi _{n}^{\left[ 0\mid \varnothing \right] \otimes \left[ \varnothing \mid
\varnothing \right] }(x;\omega ,\alpha )\propto \phi _{n}(x;\omega ,\alpha
_{1}),\quad n\in 
%TCIMACRO{\U{2124} }%
%BeginExpansion
\mathbb{Z}
%EndExpansion
.  \label{SIIO12}
\end{equation}

We can now combine the results above with the $\Gamma $ symmetries. As an example consider first the case of $\Gamma _{3}$. Combining (\ref%
{SIIO1}) and (\ref{SIIO12}) with (\ref{conjOi}) and (\ref{conjOi2}),
we have%
\begin{eqnarray}
V^{\left[ \varnothing \mid 0\right] \otimes \left[ \varnothing \mid
\varnothing \right] }(x;\omega ,\alpha ) &=&V^{\left( -1\right) \otimes
\left( \varnothing \right) }(x;\omega ,\alpha )=V(x;\omega ,\alpha )-2\left(
\log \psi _{-1}\left( x;\omega ,\alpha \right) \right) ^{\prime \prime }
\label{SIIO2} \\
&=&\Gamma _{3}\left( V(x;\omega ,\alpha )-2\left( \log \psi _{0}\left(
x;\omega ,\alpha \right) \right) ^{\prime \prime }\right) -2\omega  \notag \\
&=&\Gamma _{3}\left( V(x;\omega ,\alpha _{1})+2\omega \right) -2\omega ,
\notag 
\end{eqnarray}%
that is,

\begin{equation}
V^{\left[ \varnothing \mid 0\right] \otimes \left[ \varnothing \mid
\varnothing \right] }(x;\omega ,\alpha )=V(x;\omega ,\alpha _{-1})-2\omega .
\label{SIIO22}
\end{equation}

Moreover

\begin{eqnarray}
\psi _{n}^{\left[ \varnothing \mid 0\right] \otimes \left[ \varnothing \mid
\varnothing \right] }(x;\omega ,\alpha ) &=&\frac{W\left( \psi _{-1},\psi
_{n}\mid x;\omega ,\alpha \right) }{\psi _{-1}\left( x;\omega ,\alpha
\right) }=\Gamma _{3}\left( \frac{W\left( \psi _{0},\psi _{-\left(
n+1\right) }\mid x;\omega ,\alpha \right) }{\psi _{0}\left( x;\omega ,\alpha
\right) }\right)  \label{SIIO3} \\
&\propto &\Gamma _{3}\left( \psi _{-\left( n+1\right) -1}(x;\omega ,\alpha
_{1})\right) \propto \psi _{n+1}(x;\omega ,\alpha _{-1}),\quad n\in 
%TCIMACRO{\U{2124} }%
%BeginExpansion
\mathbb{Z}
%EndExpansion
,  \notag
\end{eqnarray}%
and

\begin{eqnarray}
\phi _{n}^{\left[ \varnothing \mid 0\right] \otimes \left[ \varnothing \mid
\varnothing \right] }(x;\omega ,\alpha ) &=&\frac{W\left( \psi _{-1},\phi
_{n}\mid x;\omega ,\alpha \right) }{\psi _{-1}\left( x;\omega ,\alpha
\right) }=\Gamma _{3}\left( \frac{W\left( \psi _{0},\phi _{-\left(
n+1\right) }\mid x;\omega ,\alpha \right) }{\psi _{0}\left( x;\omega ,\alpha
\right) }\right)  \label{SIIO31} \\
&\propto &\Gamma _{3}\left( \phi _{-\left( n+1\right) }(x;\omega ,\alpha
_{1})\right) \propto \phi _{n}(x;\omega ,\alpha _{-1}),\quad n\in 
%TCIMACRO{\U{2124} }%
%BeginExpansion
\mathbb{Z}
%EndExpansion
.  \notag
\end{eqnarray}

We can repeat the same elementary algebraic manipulation with the whole set
of $\Gamma _{i}$ symmetries and then, starting from (\ref{SIIO1}) and (%
\ref{SIIO12}), obtain in a very direct way the complete set of basic shape invariance formulas of the isotonic oscillator. For the potential, they are given by

\begin{equation}
\left\{ 
\begin{array}{l}
V^{\left[ 0\mid \varnothing \right] \otimes \left[ \varnothing \mid
\varnothing \right] }(x;\omega ,\alpha )=V(x;\omega ,\alpha
_{1})+E_{1}(\omega ) \\ 
V^{\left[ \varnothing \mid 0\right] \otimes \left[ \varnothing \mid
\varnothing \right] }(x;\omega ,\alpha )=V(x;\omega ,\alpha _{-1})-2\omega
\\ 
V^{\left[ \varnothing \mid \varnothing \right] \otimes \left[ 0\mid
\varnothing \right] }(x;\omega ,\alpha )=V(x;\omega ,\alpha _{-1}) \\ 
V^{\left[ \varnothing \mid \varnothing \right] \otimes \left[ \varnothing
\mid 0\right] }(x;\omega ,\alpha )=V(x;\omega ,\alpha _{1}),%
\end{array}%
\right.  \label{SIIOpot}
\end{equation}%
for the eigenfunctions of the extended spectrum by

\begin{equation}
\left\{ 
\begin{array}{l}
\psi _{n}^{\left[ 0\mid \varnothing \right] \otimes \left[ \varnothing \mid
\varnothing \right] }(x;\omega ,\alpha )\propto \psi _{n-1}(x;\omega ,\alpha
_{1}) \\ 
\psi _{n}^{\left[ \varnothing \mid 0\right] \otimes \left[ \varnothing \mid
\varnothing \right] }(x;\omega ,\alpha )\propto \psi _{n+1}(x;\omega ,\alpha
_{-1}) \\ 
\psi _{n}^{\left[ \varnothing \mid \varnothing \right] \otimes \left[ 0\mid
\varnothing \right] }(x;\omega ,\alpha )\propto \psi _{n}(x;\omega ,\alpha
_{-1}) \\ 
\psi _{n}^{\left[ \varnothing \mid \varnothing \right] \otimes \left[
\varnothing \mid 0\right] }(x;\omega ,\alpha )\propto \psi _{n}(x;\omega
,\alpha _{1}),%
\end{array}%
\right. \quad n\in 
%TCIMACRO{\U{2124} }%
%BeginExpansion
\mathbb{Z}
%EndExpansion
,  \label{SIIOes}
\end{equation}%
and for the eigenfunctions of the extended shadow spectrum by%
\begin{equation}
\left\{ 
\begin{array}{l}
\phi _{n}^{\left[ 0\mid \varnothing \right] \otimes \left[ \varnothing \mid
\varnothing \right] }(x;\omega ,\alpha )\propto \phi _{n}(x;\omega ,\alpha _{1})
\\ 
\phi _{n}^{\left[ \varnothing \mid 0\right] \otimes \left[ \varnothing \mid
\varnothing \right] }(x;\omega ,\alpha )\propto \phi _{n}(x;\omega ,\alpha
_{-1}) \\ 
\phi _{n}^{\left[ \varnothing \mid \varnothing \right] \otimes \left[ 0\mid
\varnothing \right] }(x;\omega ,\alpha )\propto \phi _{n-1}(x;\omega ,\alpha
_{-1}) \\ 
\phi _{n}^{\left[ \varnothing \mid \varnothing \right] \otimes \left[
\varnothing \mid 0\right] }(x;\omega ,\alpha )\propto \phi _{n+1}(x;\omega
,\alpha _{1}),%
\end{array}%
\right. \quad n\in 
%TCIMACRO{\U{2124} }%
%BeginExpansion
\mathbb{Z}
%EndExpansion
.  \label{SIIOess}
\end{equation}

Remarkably, the classical formulas for the derivative
Laguerre polynomials \cite{magnus,szego}, are a direct
transcription of the shape invariance formulas when expressed in a explicit
form. Indeed, the Darboux-Crum formula (\ref{foDBTwronsk}) allows to write (\ref{SIIOes}) as

\begin{equation}
\left\{ 
\begin{array}{l}
W\left( \psi _{0},\psi _{n}\mid x;\omega ,\alpha \right) \propto \psi
_{0}(x;\omega ,\alpha )\psi _{n-1}(x;\omega ,\alpha _{1}) \\ 
W\left( \psi _{-1},\psi _{n}\mid x;\omega ,\alpha \right) \propto \psi
_{-1}(x;\omega ,\alpha )\psi _{n+1}(x;\omega ,\alpha _{-1}) \\ 
W\left( \phi _{0},\psi _{n}\mid x;\omega ,\alpha \right) \propto \phi
_{0}(x;\omega ,\alpha )\psi _{n}(x;\omega ,\alpha _{-1}) \\ 
W\left( \phi _{0},\psi _{n}\mid x;\omega ,\alpha \right) \propto \phi
_{-1}(x;\omega ,\alpha )\psi _{n}(x;\omega ,\alpha _{1}),%
\end{array}%
\right. \quad n\in 
%TCIMACRO{\U{2124} }%
%BeginExpansion
\mathbb{Z}
%EndExpansion
.  \label{EqLagclass}
\end{equation}

Combined with the following standard properties \cite{muir} of Wronskians 

\begin{equation}
\left\{ 
\begin{array}{l}
W\left( uy_{1},...,uy_{m}\mid x\right) =u^{m}W\left( y_{1},...,y_{m}\mid
x\right) \\ 
W\left( y_{1},...,y_{m}\mid x\right) =\left( \frac{dz}{dx}\right)
^{m(m-1)/2}W\left( y_{1},...,y_{m}\mid z\right) ,%
\end{array}%
\right.  \label{wronskprop}
\end{equation}%
equation (\ref{EqLagclass}) can be rewritten as (see (\ref{spec OI}%
), (\ref{conjOi}), (\ref{shadOI}) and (\ref{conjshadOI}))

\begin{equation}
\left\{ 
\begin{array}{l}
\left( L_{n}^{\alpha }(z)\right) ^{\prime }\propto L_{n-1}^{\alpha +1}(z) \\ 
\left( z^{-\alpha }e^{z}L_{n}^{-\alpha }(-z)\right) ^{\prime } \propto
  z^{-\alpha -1}e^{z}L_{n+1}^{-\alpha -1}(-z) \\   
\left( z^{-\alpha }L_{n}^{-\alpha }(z)\right) ^{\prime }\propto z^{-\alpha
-1}L_{n}^{-\alpha -1}(z) \\ 
\left( e^{z}L_{n}^{\alpha }(-z)\right) ^{\prime }\propto e^{z}L_{n}^{\alpha
+1}(-z)%
\end{array}%
\right. .  \label{derivLag1}
\end{equation}

Taking the initial condition in $z=0$ 
\begin{equation}
L_{n}^{\alpha }(0)=\frac{\left( n+\alpha \right) \left( n+\alpha -1\right)
...\left( \alpha +1\right) }{n!},\ \left( L_{n}^{\alpha }(0)\right) ^{\prime
}=-\frac{n}{\alpha +1}L_{n}^{\alpha }(0)=-L_{n-1}^{\alpha +1}(0),
\end{equation}%
this leads to the well known formulas for the derivatives \cite{szego,magnus}%
\begin{equation}
\left\{ 
\begin{array}{l}
\left( L_{n}^{\alpha }(z)\right) ^{\prime }=-L_{n-1}^{\alpha +1}(z) \\ 
\left( z^{-\alpha }e^{z}L_{n}^{-\alpha }(-z)\right) ^{\prime }=\left(
n+1\right) z^{-\alpha -1}e^{z}L_{n+1}^{-\alpha -1}(-z) \\ 
\left( z^{-\alpha }L_{n}^{-\alpha }(z)\right) ^{\prime }=\left( n-\alpha
\right) z^{-\alpha -1}L_{n}^{-\alpha -1}(z) \\ 
\left( e^{z}L_{n}^{\alpha }(-z)\right) ^{\prime }=e^{z}L_{n}^{\alpha +1}(-z)%
\end{array}%
\right. ,  \label{derivLag12}
\end{equation}%
and more generally for the k-th derivative

\begin{equation}
\left\{ 
\begin{array}{l}
\left( L_{n}^{\alpha }(z)\right) ^{(k)}=\left( -1\right)
  ^{k}L_{n-k}^{\alpha +k}(z) \\ 
\left( z^{-\alpha }e^{z}L_{n}^{-\alpha }(-z)\right) ^{\left( k\right) }=\left(
n+1\right) _{\overline{k}}z^{-\alpha -k}e^{z}L_{n+k}^{-\alpha -k}(-z) \\ 
\left( z^{-\alpha }L_{n}^{-\alpha }(z)\right) ^{(k)}=\left( n-\alpha \right) _{\underline{k}}z^{-\alpha -k}L_{n}^{-\alpha -k}(z) \\ 
\left( e^{z}L_{n}^{\alpha }(-z)\right) ^{(k)}=e^{z}L_{n}^{\alpha +k}(-z)%
\end{array}%
\right. ,  \label{derivLag13}
\end{equation}%
where $\left( x\right) _{\overline{m}}$ and $\left( x\right) _{\underline{m}%
} $ are the usual ascending and falling factorials

\begin{equation}
\left( x\right) _{\overline{m}}=x(x+1)...(x+m-1),\ \left( x\right) _{%
\underline{m}}=x(x-1)...(x-m+1).  \label{poch}
\end{equation}

\subsection{Trigonometric Darboux-P\"oschl-Teller potential and Jacobi polynomials}

In terms of punctured universal characters, the shape invariance property (\ref{SI tDPT}) for the tDPT
potential can be written as

\begin{equation}
\left\{ 
\begin{array}{l}
V^{\left[ 0\mid \varnothing \right] \otimes \left[ \varnothing \mid
\varnothing \right] }(x;\alpha ,\beta )\equiv V^{\left( 0\right) \otimes
\left( \varnothing \right) }(x;\alpha ,\beta )=V(x;\alpha _{1},\beta
_{1})+E_{1}(\alpha ,\beta ) \\ 
\psi _{n}^{\left[ 0\mid \varnothing \right] \otimes \left[ \varnothing \mid
\varnothing \right] }(x;\alpha ,\beta )\equiv \psi _{n}^{\left( 0\right)
\otimes \left( \varnothing \right) }(x;\alpha ,\beta )\propto \psi
_{n-1}(x;\alpha _{1},\beta _{1})%
\end{array}%
\right. , \quad n\in 
%TCIMACRO{\U{2124} }%
%BeginExpansion
\mathbb{Z}
%EndExpansion
.  \label{SItDPT1}
\end{equation}%
In the particular case $n=0$, we have by (\ref{conjshadtDPT})
\begin{equation}
\psi _{0}^{\left[ \varnothing \mid 0%
\right] \otimes \left[ \varnothing \mid \varnothing \right] }(x;\alpha
,\beta )=1/\psi _{0}(x;\alpha ,\beta )=\left( 1-z\right) ^{\left( -\alpha
-1/2\right) /2}\left( 1+z\right) ^{\left( -\beta -1/2\right) /2}\propto \psi
_{-1}\left( x;\alpha _{1},\beta _{1}\right). 
\end{equation}

Observe that $\phi _{n}^{\left[ 0\mid \varnothing \right] \otimes \left[ \varnothing
\mid \varnothing \right] }(x;\alpha ,\beta )$ is a formal eigenfunction of $%
V^{\left[ 0\mid \varnothing \right] \otimes \left[ \varnothing \mid
\varnothing \right] }(x;\alpha ,\beta )$ for the formal eigenvalue $%
E_{n-\alpha }(\alpha ,\beta )$, i.e. (see above) $\phi _{n}^{\left[ 0\mid
\varnothing \right] \otimes \left[ \varnothing \mid \varnothing \right]
}(x;\alpha ,\beta )$ is a formal eigenfunction of $V(x;\alpha _{1},\beta
_{1})$ for the formal eigenvalue $E_{n-\alpha }(\alpha ,\beta )-E_{1}(\alpha
,\beta )=E_{n-\alpha _{1}}(\alpha ,\beta )$. Since it is also
quasi-polynomial we necessarily have

\begin{equation}
\phi _{n}^{\left[ 0\mid \varnothing \right] \otimes \left[ \varnothing \mid
\varnothing \right] }(x;\alpha ,\beta )\propto \phi _{n}(x;\alpha _{1},\beta
_{1}),\quad n\in 
%TCIMACRO{\U{2124} }%
%BeginExpansion
\mathbb{Z}
%EndExpansion
.  \label{SItDPT12}
\end{equation}

Proceeding as in the case of the isotonic oscillator, by combining the $\Gamma _{i}$, $i=1,2,3$
symmetries (\ref{conjtDPT}) with the results above, we obtain the following transformed potentials:

\begin{equation*}
\left\{ 
\begin{array}{l}
V^{\left[ 0\mid \varnothing \right] \otimes \left[ \varnothing \mid
\varnothing \right] }(x;\alpha ,\beta )=V(x;\alpha _{1},\beta
_{1})+E_{1}(\alpha ,\beta ) \\ 
V^{\left[ \varnothing \mid 0\right] \otimes \left[ \varnothing \mid
\varnothing \right] }(x;\alpha ,\beta )=V(x;\alpha _{-1},\beta
_{-1})+E_{-1}(\alpha ,\beta ) \\ 
V^{\left[ \varnothing \mid \varnothing \right] \otimes \left[ 0\mid
\varnothing \right] }(x;\alpha ,\beta )=V(x;\alpha _{-1},\beta _{1}) \\ 
V^{\left[ \varnothing \mid \varnothing \right] \otimes \left[ \varnothing
\mid 0\right] }(x;\alpha ,\beta )=V(x;\alpha _{1},\beta _{-1})%
\end{array}%
\right. ,
\end{equation*}%
whose eigenfunctions of the extended spectrum are

\begin{equation}
\left\{ 
\begin{array}{l}
\psi _{n}^{\left[ 0\mid \varnothing \right] \otimes \left[ \varnothing \mid
\varnothing \right] }(x;\alpha ,\beta )\propto \psi _{n-1}(x;\alpha _{1},\beta
_{1}) \\ 
\psi _{n}^{\left[ \varnothing \mid 0\right] \otimes \left[ \varnothing \mid
\varnothing \right] }(x;\alpha ,\beta )\propto \psi _{n-1}(x;\alpha _{-1},\beta
_{-1}) \\ 
\psi _{n}^{\left[ \varnothing \mid \varnothing \right] \otimes \left[ 0\mid
\varnothing \right] }(x;\alpha ,\beta )\propto \psi _{n}(x;\alpha _{-1},\beta
_{1}) \\ 
\psi _{n}^{\left[ \varnothing \mid \varnothing \right] \otimes \left[
\varnothing \mid 0\right] }(x;\alpha ,\beta )\propto \psi _{n}(x;\alpha
_{1},\beta _{-1})%
\end{array}%
\right. ,  \label{SItDPT2}
\end{equation}%
and whose eigenfunctions of the extended shadow spectrum are
\begin{equation}
\left\{ 
\begin{array}{l}
\phi _{n}^{\left[ 0\mid \varnothing \right] \otimes \left[ \varnothing \mid
\varnothing \right] }(x;\alpha ,\beta )\propto \phi _{n}(x;\alpha _{1},\beta
_{1}) \\ 
\phi _{n}^{\left[ \varnothing \mid 0\right] \otimes \left[ \varnothing \mid
\varnothing \right] }(x;\alpha ,\beta )\propto \phi _{n}(x;\alpha _{-1},\beta
_{-1}) \\ 
\phi _{n}^{\left[ \varnothing \mid \varnothing \right] \otimes \left[ 0\mid
\varnothing \right] }(x;\alpha ,\beta )\propto \phi _{n-1}(x;\alpha _{-1},\beta
_{1}) \\ 
\phi _{n}^{\left[ \varnothing \mid \varnothing \right] \otimes \left[
\varnothing \mid 0\right] }(x;\alpha ,\beta )\propto \phi _{n+1}(x;\alpha
_{1},\beta _{-1})%
\end{array}%
\right. .
\end{equation}

This exhausts the complete set of shape invariance formulas for the
quasi-polynomial eigenfunctions of the tDPT\ potential.
Combining the shape invariance formulas established above with
the Darboux-Crum formula, we can deduce in a direct way some identities for the Jacobi polynomials. In the same way as for the isotonic oscillator, we first
have

\begin{equation}
\left\{ 
\begin{array}{l}
W\left( \psi _{0},\psi _{n}\mid x;\alpha ,\beta \right) \propto \psi
_{0}(x;\alpha ,\beta )\psi _{n-1}(x;\alpha _{1},\beta _{1}) \\ 
W\left( \psi _{-1},\psi _{n}\mid x;\alpha ,\beta \right) \propto \psi
_{-1}(x;\alpha ,\beta )\psi _{n+1}(x;\alpha _{-1},\beta _{-1}) \\ 
W\left( \phi _{0},\psi _{n}\mid x;\alpha ,\beta \right) \propto \phi
_{0}(x;\alpha ,\beta )\psi _{n}(x;\alpha _{-1},\beta _{1}) \\ 
W\left( \phi _{0},\psi _{n}\mid x;\alpha ,\beta \right) \propto \phi
_{-1}(x;\alpha ,\beta )\psi _{n}(x;\alpha _{1},\beta _{-1}),%
\end{array}%
\right. \quad n\in 
%TCIMACRO{\U{2124} }%
%BeginExpansion
\mathbb{Z}
%EndExpansion
.  \label{EqJacclass}
\end{equation}

From equations (\ref{spec tDPT}), (\ref{conjtDPT}), (\ref{shadtDPT}) and (\ref%
{conjshadtDPT}), we deduce immediately

\begin{equation}
\left\{ 
\begin{array}{l}
\left( \mathit{P}_{n}^{\left( \alpha ,\beta \right) }\left( z\right) \right)
^{\prime }\propto \mathit{P}_{n-1}^{\left( \alpha +1,\beta +1\right) }\left(
z\right) \\ 
\left( \left( 1-z\right) ^{\alpha }\left( 1+z\right) ^{\beta }\mathit{P}%
_{n}^{\left( \alpha ,\beta \right) }\left( z\right) \right) ^{\prime }\propto
\left( 1-z\right) ^{\alpha -1}\left( 1+z\right) ^{\beta -1}\mathit{P}%
_{n+1}^{\left( \alpha -1,\beta -1\right) }\left( z\right) \\ 
\left( \left( 1-z\right) ^{\alpha }\mathit{P}_{n}^{\left( \alpha ,\beta
\right) }\left( z\right) \right) ^{\prime }\propto \left( 1-z\right) ^{\alpha
-1}\mathit{P}_{n}^{\left( \alpha -1,\beta +1\right) }\left( z\right) \\ 
\left( \left( 1+z\right) ^{\beta }\mathit{P}_{n}^{\left( \alpha ,\beta
\right) }\left( z\right) \right) ^{\prime }\propto \left( 1+z\right) ^{\beta -1}%
\mathit{P}_{n}^{\left( \alpha +1,\beta -1\right) }\left( z\right)%
\end{array}%
\right. ,  \label{derivJac10}
\end{equation}

Evaluating the above expressions in $z=1$, we have
\begin{equation}
\mathit{P}_{n}^{\left( \alpha ,\beta \right) }\left( 1\right) =\frac{\left(
n+\alpha \right) \left( n+\alpha -1\right) ...\left( \alpha +1\right) }{n!},%
\mathit{P}_{n}^{\left( \alpha ,\beta \right) \prime }\left( 1\right) =\frac{%
n+\alpha +\beta +1}{2}\mathit{P}_{n-1}^{\left( \alpha +1,\beta +1\right)
}\left( 1\right) .
\end{equation}%
and making use of the reflection formula for Jacobi polynomials \cite%
{szego,magnus} 
\[\mathit{P}%
_{n}^{\left( \beta ,\alpha \right) }\left( -z\right) =\left( -1\right) ^{n}%
\mathit{P}_{n}^{\left( \alpha ,\beta \right) }\left( z\right), \]
 we obtain

\begin{equation}
\left\{ 
\begin{array}{l}
\mathit{P}_{n}^{\left( \alpha ,\beta \right) \prime }\left( z\right) =\frac{%
n+\alpha +\beta +1}{2}\mathit{P}_{n-1}^{\left( \alpha +1,\beta +1\right)
}\left( z\right) \\ 
\left( \left( 1-z\right) ^{\alpha }\left( 1+z\right) ^{\beta }\mathit{P}%
_{n}^{\left( \alpha ,\beta \right) }\left( z\right) \right) ^{\prime
}=-2\left( n+1\right) \left( 1-z\right) ^{\alpha -1}\left( 1+z\right)
^{\beta -1}\mathit{P}_{n+1}^{\left( \alpha -1,\beta -1\right) }\left(
z\right) \\ 
\left( \left( 1-z\right) ^{\alpha }\mathit{P}_{n}^{\left( \alpha ,\beta
\right) }\left( z\right) \right) ^{\prime }=-\left( n+\alpha \right) \left(
1-z\right) ^{\alpha -1}\mathit{P}_{n}^{\left( \alpha -1,\beta +1\right)
}\left( z\right) \\ 
\left( \left( 1+z\right) ^{\beta }\mathit{P}_{n}^{\left( \alpha ,\beta
\right) }\left( z\right) \right) ^{\prime }=\left( n+\beta \right) \left(
1+z\right) ^{\beta -1}\mathit{P}_{n}^{\left( \alpha +1,\beta -1\right)
}\left( z\right)%
\end{array}%
\right. .  \label{derivJac101}
\end{equation}

\bigskip More generally, for the k-th derivative we have (\ref{poch})

\begin{equation}
\left\{ 
\begin{array}{l}
\left( \mathit{P}_{n}^{\left( \alpha ,\beta \right) }\left( z\right) \right)
^{\left( k\right) }=\frac{1}{2^{k}}\left( n+\alpha +\beta +1\right) _{k}%
\mathit{P}_{n-k}^{\left( \alpha +k,\beta +k\right) }\left( z\right) \\ 
\left( \left( 1-z\right) ^{\alpha }\left( 1+z\right) ^{\beta }\mathit{P}%
_{n}^{\left( \alpha ,\beta \right) }\left( z\right) \right) ^{\left(
k\right) }=\left( -2\right) ^{k}\left( n+1\right) _{\overline{k}}\left(
1-z\right) ^{\alpha -k}\left( 1+z\right) ^{\beta -k}\mathit{P}_{n+k}^{\left(
\alpha -k,\beta -k\right) }\left( z\right) \\ 
\left( \left( 1-z\right) ^{\alpha }\mathit{P}_{n}^{\left( \alpha ,\beta
\right) }\left( z\right) \right) ^{\left( k\right) }=\left( -1\right)
^{k}\left( n+\alpha \right) _{\underline{k}}\left( 1-z\right) ^{\alpha -k}%
\mathit{P}_{n}^{\left( \alpha -k,\beta +k\right) }\left( z\right) \\ 
\left( \left( 1+z\right) ^{\beta }\mathit{P}_{n}^{\left( \alpha ,\beta
\right) }\left( z\right) \right) ^{\left( k\right) }=\left( n+\beta \right)
_{\underline{k}}\left( 1+z\right) ^{\beta -k}\mathit{P}_{n}^{\left( \alpha
+k,\beta -k\right) }\left( z\right)%
\end{array}%
\right. .  \label{derivJac102}
\end{equation}

We see thus that the classical formulas for the derivatives of Jacobi polynomials 
\cite{magnus,szego}, are just an explicit transcription of
the complete set of shape invariance formulas for the eigenfunctions of the
tDPT potential.

%%%%%%%%%%%%%%%%%%%%%%%%%%%%%%%%%%%%%%%%%%%%%%%%%%%%%%%%%%%%%%%%%%%
\section{Equivalence of Laguerre pseudowronskians}\label{sec:L}
%%%%%%%%%%%%%%%%%%%%%%%%%%%%%%%%%%%%%%%%%%%%%%%%%%%%%%%%%%%%%%%%%%%

In this section we shall derive many identities among determinants whose entries are associated Laguerre polynomials, by exploiting the notion of equivalence of punctured universal characters explained in Section \ref{sec:PUC}.

Consider a mixed chain associated to the bi-tuple $N_{m}\otimes L_{r}=\left(
n_{1},...,n_{k},\overline{n}_{\overline{k}},...,\overline{n}_{1}\right)
\otimes \left( l_{1},...,l_{r}\right) $.

We consider first the case where $n_{i}>n_{i+1}>0$ and $-1\geq \overline{n}%
_{i+1}>\overline{n}_{i}$. We add to the first Maya diagram the index $0$ (ie
we add to the preceding chain a Darboux transformation associated to the seed function $\psi
_{0} $) obtaining the bi-tuple $\left( 0,N_{m}\right) \otimes L_{r}$.

Considering $V^{\left( 0,N_{m}\right) \otimes L_{r}}$ as an extension of $%
V^{\left( 0\right) \otimes \varnothing }$, the Crum formula (\ref{etats n3}%
) applied to the successive seed functions along the chain gives, after
telescopic cancellation:

\begin{eqnarray*}
W\left( \psi _{n_{1}}^{\left( 0\right) \otimes \varnothing },...,\psi _{%
\overline{n}_{1}}^{\left( 0\right) \otimes \varnothing },\phi
_{l_{1}}^{\left( 0\right) \otimes \varnothing },...,\phi _{l_{r}}^{\left(
0\right) \otimes \varnothing }\mid x;\omega ,\alpha \right) &=&\psi
_{n_{1}}^{\left( 0\right) \otimes \varnothing }\left( x;\omega ,\alpha
\right) ...\psi _{\overline{n}_{1}}^{\left( 0,N_{m}\right) \otimes
\varnothing }\left( x;\omega ,\alpha \right)  \label{a1} \\
&&\times \phi _{l_{1}}^{\left( 0,N_{m}\right) \otimes \varnothing }\left(
x;\omega ,\alpha \right) ...\phi _{l_{r}}^{\left( 0,N_{m}\right) \otimes
L_{r}}\left( x;\omega ,\alpha \right).  \notag
\end{eqnarray*}

On the other hand, if we consider $V^{\left( 0,N_{m}\right) \otimes L_{r}}$ as
an extension of $V$ , the same use of Crum's formulas (\ref{etats n3})
leads after telescopic cancellation to

\begin{eqnarray*}
W\left( \psi _{0},\psi _{n_{1}},...,\psi _{\overline{n}_{1}},\phi
_{l_{1}},...,\phi _{l_{r}}\mid x;\omega ,\alpha \right) &=&\psi _{0}\left(
x;\omega ,\alpha \right) \psi _{n_{1}}^{\left( 0\right) \otimes \varnothing
}\left( x;\omega ,\alpha \right) ...\psi _{\overline{n}_{1}}^{\left(
0,N_{m}\right) \otimes \varnothing }\left( x;\omega ,\alpha \right)
\label{a2} \\
&&\times \phi _{l_{1}}^{\left( 0,N_{m}\right) \otimes \varnothing }\left(
x;\omega ,\alpha \right) ...\phi _{l_{r}}^{\left( 0,N_{m}\right) \otimes
L_{r}}\left( x;\omega ,\alpha \right) .  \notag
\end{eqnarray*}

Comparing the last two identities and making use of shape invariance
(\ref{SIIOes}) and Eq(\ref{SIIOess}), we obtain immediately

\begin{equation}
W^{\left( 0,N_{m}\right) \otimes L_{r}}\left( x;\omega ,\alpha \right) \propto
\psi _{0}\left( x;\omega ,\alpha \right) W^{\left( N_{m}-1\right) \otimes
L_{r}}\left( x;\omega ,\alpha _{1}\right) ,  \label{eqwroniso1}
\end{equation}%
or, using (\ref{psi-1}) and redefining $\alpha \rightarrow \alpha _{-1}$

\begin{equation}
W^{\left( N_{m}-1\right) \otimes L_{r}}\left( x;\omega ,\alpha \right) \propto
\psi _{-1}\left( x;\omega ,\alpha \right) W^{\left( 0,N_{m}\right) \otimes
\varnothing }\left( x;\omega ,\alpha _{-1}\right).  \label{eqwroniso12}
\end{equation}

Suppose $\overline{n}_{\overline{k}}<-1$. By shifting the indices $%
n_{i}\rightarrow n_{i}+1,$ $\overline{n}_{i}\rightarrow \overline{n}_{i}+1$,
(\ref{eqwroniso12}) can be rewritten as (see (\ref{shiftn3})) 
\begin{equation}
W^{N_{m}\otimes L_{r}}\left( x;\omega ,\alpha \right) \propto \psi _{-1}\left(
x;\omega ,\alpha \right) W^{\left( N_{m}\right) _{1}\otimes L_{r}}\left(
x;\omega ,\alpha _{-1}\right) .  \label{eqiso1}
\end{equation}

Suppose now that in $N_{m}$ we have $n_{k}=0$ (in which case $n_{k-1}>0$).
(\ref{eqwroniso1}) is still valid provided we change $k$ into $k-1$. Then

\begin{equation}
W^{\left( 0,N_{m}\right) \otimes L_{r}}\left( x;\omega ,\alpha \right) \propto
\psi _{0}\left( x;\omega ,\alpha \right) W^{\left( N_{m}-1\right) \otimes
L_{r}}\left( x;\omega ,\alpha _{1}\right) ,
\end{equation}%
i.e. with the shifting rule (\ref{shiftn2})

\begin{equation}
W^{N_{m}\otimes L_{r}}\left( x;\omega ,\alpha \right) \propto \psi _{0}\left(
x;\omega ,\alpha \right) W^{_{1}\left( N_{m}\right) \otimes L_{r}}\left(
x;\omega ,\alpha _{1}\right) .
\end{equation}

Consider in a second step, a mixed chain associated to the bi-tuple of
relative integers $N_{m}\otimes L_{r}$ to which we add to the Maya diagram $%
N_{m}$ the index $(-1)$ (i.e. we add a DT associated to the seed function $%
\psi _{-1}$ to the corresponding chain) obtaining $\left( -1,N_{m}\right)
\otimes L_{r}$. Proceeding as above, using the Darboux-Crum and Crum
formulas as well as telescopic cancellation and the shape invariance
properties (\ref{SIIOes}) and (\ref{SIIOess}), we obtain

\begin{equation}
W^{\left( -1,N_{m}\right) \otimes L_{r}}\left( x;\omega ,\alpha \right) \propto
\psi _{-1}\left( x;\omega ,\alpha \right) W^{\left( N_{m}+1\right) \otimes
L_{r}}\left( x;\omega ,\alpha _{-1}\right) .  \label{eqwroniso2}
\end{equation}

If $n_{k}>0$, shifting the indices ($n_{i}\rightarrow n_{i}-1,$ $\overline{n}%
_{i}\rightarrow \overline{n}_{i}-1$) and using (\ref{psi-1}) we have
equivalently 
\begin{equation}
W^{N_{m}\otimes L_{r}}\left( x;\omega ,\alpha \right) \propto \psi _{0}\left(
x;\omega ,\alpha \right) W^{_{1}\left( N_{m}\right) \otimes L_{r}}\left(
x;\omega ,\alpha _{1}\right) .  \label{eqiso2}
\end{equation}

If $\overline{n}_{\overline{k}}=-1$ (in which case $\overline{n}_{\overline{k%
}-1}<-1$), (\ref{eqwroniso2}) gives

\begin{equation}
W^{\left( -1,N_{m}\right) \otimes L_{r}}\left( x;\omega ,\alpha \right) \propto
\psi _{-1}\left( x;\omega ,\alpha \right) W^{\left( N_{m}+1\right) \otimes
L_{r}}\left( x;\omega ,\alpha _{-1}\right),
\end{equation}%
i.e. (see Eq(\ref{shiftn2}))

\begin{equation}
W^{N_{m}\otimes L_{r}}\left( x;\omega ,\alpha \right) \propto \psi _{-1}\left(
x;\omega ,\alpha \right) W^{\left( N_{m}\right) _{1}\otimes L_{r}}\left(
x;\omega ,\alpha _{-1}\right) .
\end{equation}

Summarizing, we have for every bi-tuple $N_{m}\otimes L_{r}$

\begin{equation}
\left\{ 
\begin{array}{l}
W^{N_{m}\otimes L_{r}}\left( x;\omega ,\alpha \right) \propto \psi _{0}\left(
x;\omega ,\alpha \right) W^{_{1}\left( N_{m}\right) \otimes L_{r}}\left(
x;\omega ,\alpha _{1}\right) \\ 
W^{N_{m}\otimes L_{r}}\left( x;\omega ,\alpha \right) \propto \psi _{-1}\left(
x;\omega ,\alpha \right) W^{\left( N_{m}\right) _{1}\otimes L_{r}}\left(
x;\omega ,\alpha _{-1}\right) .%
\end{array}%
\right.  \label{eq1step1}
\end{equation}

We can repeat exactly the same lines by acting on the Maya diagram $L_{r}$
rather than on $N_{m}$. We obtain in the same way

\begin{equation}
\left\{ 
\begin{array}{l}
W^{N_{m}\otimes L_{r}}\left( x;\omega ,\alpha \right) \propto \phi _{0}\left(
x;\omega ,\alpha \right) W^{N_{m}\otimes \ _{1}\left( L_{r}\right) }\left(
x;\omega ,\alpha _{1}\right) \\ 
W^{N_{m}\otimes L_{r}}\left( x;\omega ,\alpha \right) \propto \phi _{-1}\left(
x;\omega ,\alpha \right) W^{N_{m}\otimes \left( L_{r}\right) _{1}}\left(
x;\omega ,\alpha _{1}\right) .%
\end{array}%
\right.  \label{eq1step2}
\end{equation}

Note that the identities (\ref{EqLagclass}) and the 
formulas for the derivatives of Laguerre polynomials (\ref{EqLagclass}) constitute in fact the lowest order case of (\ref{eq1step1}) and (\ref{eq1step2}).

If we look at the consequences of the relations (\ref{eq1step1}) and (%
\ref{eq1step2}) at the level of the extended potentials, we obtain first,
using (\ref{SIIOpot})

\begin{eqnarray}
V^{N_{m}\otimes L_{r}}(x;\omega ,\alpha ) &=&V(x;\omega ,\alpha )-2\left(
\log \left( \psi _{0}\left( x;\omega ,\alpha \right) W^{_{1}\left(
N_{m}\right) \otimes L_{r}}\left( x;\omega ,\alpha _{1}\right) \right)
\right) ^{\prime \prime } \\
&=&V^{\left( 0\right) \otimes \varnothing }(x;\omega ,\alpha )-2\left( \log
W^{_{1}\left( N_{m}\right) \otimes L_{r}}\left( x;\omega ,\alpha _{1}\right)
\right) ^{\prime \prime }  \notag \\
&=&V(x;\omega ,\alpha _{1})+E_{1}(\omega )-2\left( \log W^{_{1}\left(
N_{m}\right) \otimes L_{r}}\left( x;\omega ,\alpha _{1}\right) \right)
^{\prime \prime }  \notag \\
&=&V^{_{1}\left( N_{m}\right) \otimes L_{r}}(x;\omega ,\alpha
_{1})+E_{1}(\omega ).  \notag
\end{eqnarray}

More generally, we have

\begin{equation}
\left\{ 
\begin{array}{l}
V^{N_{m}\otimes L_{r}}(x;\omega ,\alpha )=V^{_{1}\left( N_{m}\right) \otimes
L_{r}}(x;\omega ,\alpha _{1})+E_{1}(\omega ) \\ 
V^{N_{m}\otimes L_{r}}(x;\omega ,\alpha )=V^{\left( N_{m}\right) _{1}\otimes
L_{r}}(x;\omega ,\alpha _{-1})+E_{-1}(\omega ) \\ 
V^{N_{m}\otimes L_{r}}(x;\omega ,\alpha )=V^{N_{m}\otimes \ _{1}\left(
L_{r}\right) }(x;\omega ,\alpha _{1}) \\ 
V^{N_{m}\otimes L_{r}}(x;\omega ,\alpha )=V^{N_{m}\otimes \left(
L_{r}\right) _{1}}(x;\omega ,\alpha _{-1}).%
\end{array}%
\right.  \label{eq1steppot}
\end{equation}

In terms of punctured universal characters the preceding results (\ref{eq1step1}), (\ref{eq1step2}) and (\ref{eq1steppot}) can be rewritten as

\begin{equation}
\left\{ 
\begin{array}{l}
W^{\left[ d\mid \overline{d}\right] \otimes \left[ \delta \mid \overline{%
\delta }\right] }\left( x;\omega ,\alpha \right) \propto \psi _{0}\left(
x;\omega ,\alpha \right) W^{_{1}\left[ d\mid \overline{d}\right] \otimes %
\left[ \delta \mid \overline{\delta }\right] }\left( x;\omega ,\alpha
_{1}\right) \\ 
W^{\left[ d\mid \overline{d}\right] \otimes \left[ \delta \mid \overline{%
\delta }\right] }\left( x;\omega ,\alpha \right) \propto \psi _{-1}\left(
x;\omega ,\alpha \right) W^{\left[ d\mid \overline{d}\right] _{1}\otimes %
\left[ \delta \mid \overline{\delta }\right] }\left( x;\omega ,\alpha
_{-1}\right) \\ 
W^{\left[ d\mid \overline{d}\right] \otimes \left[ \delta \mid \overline{%
\delta }\right] }\left( x;\omega ,\alpha \right) \propto \phi _{0}\left(
x;\omega ,\alpha \right) W^{\left[ d\mid \overline{d}\right] \otimes \ _{1}%
\left[ \delta \mid \overline{\delta }\right] }\left( x;\omega ,\alpha
_{-1}\right) \\ 
W^{\left[ d\mid \overline{d}\right] \otimes \left[ \delta \mid \overline{%
\delta }\right] }\left( x;\omega ,\alpha \right) \propto \phi _{-1}\left(
x;\omega ,\alpha \right) W^{\left[ d\mid \overline{d}\right] \otimes \left[
\delta \mid \overline{\delta }\right] _{1}}\left( x;\omega ,\alpha
_{1}\right)%
\end{array}%
\right.
\end{equation}%
and

\begin{equation}
\left\{ 
\begin{array}{l}
V^{\left[ d\mid \overline{d}\right] \otimes \left[ \delta \mid \overline{%
\delta }\right] }(x;\omega ,\alpha )=V^{_{1}\left[ d\mid \overline{d}\right]
\otimes \left[ \delta \mid \overline{\delta }\right] }(x;\omega ,\alpha
_{1})+E_{1}(\omega ) \\ 
V^{\left[ d\mid \overline{d}\right] \otimes \left[ \delta \mid \overline{%
\delta }\right] }(x;\omega ,\alpha )=V^{\left[ d\mid \overline{d}\right]
_{1}\otimes \left[ \delta \mid \overline{\delta }\right] }(x;\omega ,\alpha
_{-1})+E_{-1}(\omega ) \\ 
V^{\left[ d\mid \overline{d}\right] \otimes \left[ \delta \mid \overline{%
\delta }\right] }(x;\omega ,\alpha )=V^{\left[ d\mid \overline{d}\right]
\otimes \ _{1}\left[ \delta \mid \overline{\delta }\right] }(x;\omega
,\alpha _{-1}) \\ 
V^{\left[ d\mid \overline{d}\right] \otimes \left[ \delta \mid \overline{%
\delta }\right] }(x;\omega ,\alpha )=V^{\left[ d\mid \overline{d}\right]
\otimes \left[ \delta \mid \overline{\delta }\right] _{1}}(x;\omega ,\alpha
_{1}),%
\end{array}%
\right.
\end{equation}%
where $\left[ d\mid \overline{d}\right] \otimes \left[ \delta \mid \overline{%
\delta }\right] $ ($\left[ d\mid \overline{d}\right] =\left[
d_{1},...,d_{k}\mid \overline{d}_{1},...,\overline{d}_{\overline{k}}\right] $
with $k+\overline{k}=m$, and $\left[ \delta \mid \overline{\delta }\right] =%
\left[ \delta _{1},...,\delta _{\kappa }\mid \overline{\delta }_{1},...,%
\overline{\delta }_{\overline{\kappa }}\right] $ with $\kappa +\overline{%
\kappa }=r$) is the punctured universal character corresponding to the bi-tuple $N_{m}\otimes L_{r}$.

By repeating the procedure we can push the origins in both punctured Young
diagrams in the upper right corner.

Indeed, by\ doing $-\overline{n}_{1}=\overline{d}_{1}+\overline{k}$
successive shifts in the right-ascending direction in the first Young
diagram and $\overline{l}_{1}+1=\overline{\delta }_{1}+\overline{\kappa }$
shifts in the same direction in the second one, we obtain

\begin{equation}
W^{\left[ d\mid \overline{d}\right] \otimes \left[ \delta \mid \overline{%
\delta }\right] }\left( x;\omega ,\alpha \right) \propto \prod_{i=0}^{\overline{d%
}_{1}+\overline{k}-1}\psi _{-1}\left( x;\omega ,\alpha _{-i}\right)
\prod_{j=0}^{\overline{\delta }_{1}+\overline{\kappa }-1}\phi _{-1}\left(
x;\omega ,\alpha _{j-\left( \overline{d}_{1}+\overline{k}\right) }\right) W^{%
\left[ \lambda \mid \varnothing \right] \otimes \left[ \mu \mid \varnothing %
\right] }\left( x;\omega ,\alpha _{\overline{\delta }_{1}+\overline{\kappa }-%
\overline{d}_{1}-\overline{k}}\right) .  \label{EqIO2}
\end{equation}

Equivalently we can choose to push the origin on each punctured Young
diagram in the lower left corner by\ doing $n_{1}+1=d_{1}+k$ successive
shifts in the left-descending direction in the first Young diagram and $%
l_{1}+1=\delta _{1}+\kappa $ shifts in the same direction in the second one,
obtain then

\begin{equation}
W^{\left[ d\mid \overline{d}\right] \otimes \left[ \delta \mid \overline{%
\delta }\right] }\left( x;\omega ,\alpha \right) \propto
\prod_{i=0}^{d_{1}+k-1}\psi _{0}\left( x;\omega ,\alpha _{i}\right)
\prod_{j=0}^{\delta _{1}+\kappa -1}\phi _{0}\left( x;\omega ,\alpha
_{d_{1}+k-j}\right) W^{\left[ \varnothing \mid \overline{\lambda }\right]
\otimes \left[ \varnothing \mid \overline{\mu }\right] }\left( x;\omega
,\alpha _{d_{1}+k-(\delta _{1}+\kappa )}\right) .  \label{EqIO}
\end{equation}
In the same manner we arrive at

\begin{equation}
\left\{ 
\begin{array}{l}
W^{\left[ d\mid \overline{d}\right] \otimes \left[ \delta \mid \overline{%
\delta }\right] }\left( x;\omega ,\alpha \right) \propto \prod_{i=0}^{\overline{d%
}_{1}+\overline{k}-1}\psi _{-1}\left( x;\omega ,\alpha _{-i}\right)
\prod_{j=0}^{\delta _{1}+\kappa -1}\phi _{0}\left( x;\omega ,\alpha _{-\left(
d_{1}+k+j\right) }\right) W^{\left[ \lambda \mid \varnothing \right] \otimes %
\left[ \varnothing \mid \overline{\mu }\right] }\left( x;\omega ,\alpha _{-(%
\overline{d}_{1}+\overline{k}+\delta _{1}+\kappa )}\right) \\ 
W^{\left[ d\mid \overline{d}\right] \otimes \left[ \delta \mid \overline{%
\delta }\right] }\left( x;\omega ,\alpha \right) \propto
\prod_{i=0}^{d_{1}+k-1}\psi _{0}\left( x;\omega ,\alpha _{i}\right)
\prod_{j=0}^{\overline{\delta }_{1}+\overline{\kappa }-1}\phi _{-1}\left(
x;\omega ,\alpha _{d_{1}+k+j}\right) W^{\left[ \varnothing \mid \overline{%
\lambda }\right] \otimes \left[ \varnothing \mid \overline{\mu }\right]
}\left( x;\omega ,\alpha _{d_{1}+k+\overline{\delta }_{1}+\overline{\kappa }%
}\right)%
\end{array}%
\right. .
\end{equation}
For the potentials, this gives

\begin{equation}
\left\{ 
\begin{array}{l}
V^{\left[ d\mid \overline{d}\right] \otimes \left[ \delta \mid \overline{%
\delta }\right] }\left( x;\omega ,\alpha \right) =V^{\left[ \lambda \mid
\varnothing \right] \otimes \left[ \mu \mid \varnothing \right] }\left(
x;\omega ,\alpha _{\overline{\delta }_{1}+\overline{\kappa }-\overline{d}%
_{1}-\overline{k}}\right) +E_{-\left( \overline{d}_{1}+\overline{k}\right)
}(\omega ) \\ 
V^{\left[ d\mid \overline{d}\right] \otimes \left[ \delta \mid \overline{%
\delta }\right] }\left( x;\omega ,\alpha \right) =V^{\left[ \varnothing \mid 
\overline{\lambda }\right] \otimes \left[ \varnothing \mid \overline{\mu }%
\right] }(x;\omega ,\alpha _{d_{1}+k-\delta _{1}-\kappa
})+E_{d_{1}+k}(\omega ).%
\end{array}%
\right.  \label{Equivpotiso}
\end{equation}

These formulas constitute the general transcription of the shape invariance
property of the isotonic potential at the level of its rationally extended
descendants.

In the following we choose as canonical representative of a given
shape class the punctured universal character $\left[ \lambda \mid \varnothing \right] \otimes \left[
\mu \mid \varnothing \right] $ and then express all the extensions of this
class in terms of $W^{\left[ \lambda \mid \varnothing \right] \otimes \left[
\mu \mid \varnothing \right] }$. Using (\ref{conjshadOI}) and (\ref%
{conjOi}), we obtain a more explicit form for (\ref{EqIO2})

\begin{eqnarray}
W^{\left[ d\mid \overline{d}\right] \otimes \left[ \delta \mid \overline{%
\delta }\right] }\left( x;\omega ,\alpha \right) &\propto &z^{-\alpha (%
\overline{d}_{1}+\overline{k}-\overline{\delta }_{1}-\overline{\kappa }%
)/2}\times e^{z(\overline{d}_{1}+\overline{k}+\overline{\delta }_{1}+%
\overline{\kappa })/2}  \label{EqIO3} \\
&&\times z^{(\overline{d}_{1}+\overline{k}-\overline{\delta }_{1}-\overline{%
\kappa })^{2}/4}\times W^{\left[ \lambda \mid \varnothing \right] \otimes %
\left[ \mu \mid \varnothing \right] }\left( x;\omega ,\alpha _{-(\overline{d}%
_{1}+\overline{k})+\overline{\delta }_{1}+\overline{\kappa }}\right) ,
\notag
\end{eqnarray}%
where $z=\omega x^{2}/2$.

From (\ref{conjOi}), (\ref{shadOI}) and (\ref{conjshadOI}) we have
\begin{equation}
\left\{ 
\begin{array}{l}
\psi _{-\left( n+1\right) }\left( x;\omega ,\alpha \right) =z^{-\alpha
}e^{z}\psi _{0}\left( x;\omega ,\alpha \right) L_{n}^{-\alpha }(-z) \\ 
\phi _{n}\left( x;\omega ,\alpha \right) =z^{-\alpha }\psi _{0}\left(
x;\omega ,\alpha \right) L_{n}^{-\alpha }(z) \\ 
\phi _{-\left( n+1\right) }\left( x;\omega ,\alpha \right) =e^{z}\psi
_{0}\left( x;\omega ,\alpha \right) L_{n}^{\alpha }(-z)\,.%
\end{array}%
\right. 
\end{equation}
Then using (\ref{EqIO}), (\ref{EqIO2}) and (\ref{derivLag13}), and noting that
 $m=k+\overline{k},r=\kappa +\overline{\kappa }$, we deduce that

\begin{eqnarray}
W^{\left[ d\mid \overline{d}\right] \otimes \left[ \delta \mid \overline{%
\delta }\right] }\left( x;\omega ,\alpha \right)  &\propto &z^{\alpha \left( k-%
\overline{k}-\kappa +\overline{\kappa }\right) /2}e^{-z\left( k-\overline{k}%
+\kappa -\overline{\kappa }\right) /2}  \label{W1iso} \\
&&\times z^{\left( k-\overline{k}-\kappa +\overline{\kappa }\right)
^{2}/4}z^{-\left( \overline{k}+\kappa \right) \left( \overline{k}+\kappa
-1\right) }\times L^{\left[ d\mid \overline{d}\right] \otimes \left[ \delta
\mid \overline{\delta }\right] }\left( z;\alpha \right) ,  \notag
\end{eqnarray}%
where $L^{\left[ d\mid \overline{d}\right] \otimes \left[ \delta \mid 
\overline{\delta }\right] }\left( z;\alpha \right) $ is a polynomial which we shall denote as a \textit{Laguerre pseudowronskian}. A Laguerre pseudowronskian is defined by the following determinantal expression

\begin{equation}
L^{\left[ d\mid \overline{d}\right] \otimes \left[ \delta \mid \overline{%
\delta }\right] }\left( z;\alpha \right) =\det \left( \left\{ 
\overrightarrow{\mathcal{L}}^{\left( j\right) }\left( z,\alpha \right)
\right\} _{1\leq j\leq k},\left\{ \overrightarrow{\overline{\mathcal{L}}}%
^{\left( j\right) }\left( z,\alpha \right) \right\} _{1\leq j\leq \overline{k%
}},\left\{ \overrightarrow{\Lambda }^{\left( j\right) }\left( z,\alpha
\right) \right\} _{1\leq j\leq \kappa },\left\{ \overrightarrow{\overline{%
\Lambda }}^{\left( j\right) }\left( z,\alpha \right) \right\} _{1\leq
j\leq  \overline{\kappa }}\right) ,  \label{XLP1}
\end{equation}%
where

\begin{equation}
\left\{ 
\begin{array}{l}
  \mathcal{L}_{i}^{\left( j\right) }\left( z,\alpha \right) =\left( -1\right)
  ^{i-1}L_{n_{j}-i+1}^{\alpha +i-1}\left( z\right) \\ 
  \overline{\mathcal{L}}_{i}^{\left( j\right) }\left( z,\alpha \right) =\left(
  -\overline{n}_{j}\right) _{\overline{i-1}}\times z^{m+r-i}L_{-\overline{n}%
  _{j}+i-2}^{-\alpha -i+1}\left( -z\right) \\ 
  \Lambda _{i}^{\left( j\right) }\left( z,\alpha \right) =\left( l_{j}+\alpha
  \right) _{\underline{i-1}}\times z^{m+r-i}L_{l_{j}}^{-\alpha -i+1}\left(
  z\right) \\ 
  \overline{\Lambda }_{i}^{\left( j\right) }\left( z,\alpha \right) =L_{-\overline{l}_{j}-1}^{\alpha +i-1}\left( -z\right)%
\end{array}%
\right. ,\ i=1,...,m+r,  \label{XLP2}
\end{equation}%
with $n_{j}=d_{j}+k-j\geq 0$, $l_{j}=\delta _{j}+\kappa -j\geq 0$, $-%
\overline{n}_{j}-1=\overline{d}_{j}+\overline{k}-j\geq 0$ and $-\overline{l}%
_{j}-1=\overline{\delta }_{j}+\overline{\kappa }-j\geq 0$.

If $\left[ \lambda \mid \varnothing \right] $ is in the same shape class as $%
\left[ d\mid \overline{d}\right] =\left[ d_{1},...,d_{k}\mid \overline{d}%
_{1},...,\overline{d}_{\overline{k}}\right] $, we can use (\ref{W1iso})
replacing $\overline{k}$ and $\overline{\kappa }$ by $0$ and $k$ and $\kappa 
$ by $l(\lambda )=\overline{d}_{1}+k$ and $l(\mu )$ $=\overline{\delta }%
_{1}+\kappa $ respectively. It gives

\begin{eqnarray}
W^{\left[ \lambda \mid \varnothing \right] \otimes \left[ \mu \mid
\varnothing \right] }\left( x;\omega ,\alpha \right) &\propto &z^{\alpha \left( 
\overline{d}_{1}+k-\overline{\delta }_{1}-\kappa \right) /2}\times
e^{-z\left( \overline{d}_{1}+k+\overline{\delta }_{1}+\kappa \right) /2}
\label{W2iso} \\
&&\times z^{\left( \overline{d}_{1}+k-\overline{\delta }_{1}-\kappa \right)
^{2}/4}\times z^{-(\overline{\delta }_{1}+\kappa )\left( \overline{\delta }%
_{1}+\kappa -1\right) }\times L^{\left[ \lambda \mid \varnothing \right]
\otimes \left[ \mu \mid \varnothing \right] }\left( z;\alpha \right) . 
\notag
\end{eqnarray}

The equivalence relation (\ref{EqIO3}) can then be rewritten as
\begin{eqnarray*}
W^{\left[ d\mid \overline{d}\right] \otimes \left[ \delta \mid \overline{%
\delta }\right] }\left( x;\omega ,\alpha \right) &\propto &z^{\alpha (k-%
\overline{k}-\kappa +\overline{\kappa })/2}\times e^{-z\left( k-\overline{k}%
+\kappa -\overline{\kappa }\right) /2} \\
&&\times z^{(k-\overline{k}-\kappa +\overline{\kappa })^{2}/4}\times z^{-(%
\overline{\delta }_{1}+\kappa )\left( \overline{\delta }_{1}+\kappa
-1\right) }\times L^{\left[ \lambda \mid \varnothing \right] \otimes \left[
\mu \mid \varnothing \right] }\left( z;\alpha _{-\left( \overline{d}_{1}+%
\overline{k}\right) +\overline{\delta }_{1}+\overline{\kappa }}\right) , 
\notag
\end{eqnarray*}%
and finally, using Eq(\ref{W1iso}) we arrive at

\begin{equation}
L^{\left[ d\mid \overline{d}\right] \otimes \left[ \delta \mid \overline{%
\delta }\right] }\left( z;\alpha \right) \propto z^{\left( \overline{k}+\kappa
\right) \left( \overline{k}+\kappa -1\right) -(\overline{\delta }_{1}+\kappa
)\left( \overline{\delta }_{1}+\kappa -1\right) }\times L^{\left[ \lambda
\mid \varnothing \right] \otimes \left[ \mu \mid \varnothing \right] }\left(
z;\alpha _{-\left( \overline{d}_{1}+\overline{k}\right) +\overline{\delta }%
_{1}+\overline{\kappa }}\right) ,  \label{EquivXLP}
\end{equation}%
which is the \textit{equivalence relation for Laguerre pseudowronskians}.

\begin{example}

  As an example of an equivalence relation, consider the case where
  $\left[ d\mid \overline{d}\right] =\left[ 2,1\mid 1\right] $,
  $N_{3}=(3,1,-2)$ and $\left[ \delta \mid \overline{%
      \delta }\right] =\left[ 3\mid 2\right] $,
  $L_{2}=(3,-3)$. 
Setting 
\begin{align}
  \cL_{i,j} &= (-1)^i L^{(\alpha+i)}_j(z) \\
  \bcL_{i,j} &=  \ffac{j+i}{i} L^{(-\alpha-i)}_j(-z) \\
  \Lambda_{i,j} &=   \ffac{j-\alpha}{i} L^{(-\alpha-i)}_j(z)\\
  \bLambda_{i,j} &=  L^{(\alpha+i)}_j(-z),
\end{align}
 the pseudo-Wronskians in question have the explicit
determinantal form
\begin{gather}
  L^{\left[ 2,1\mid 1\right] \otimes \left[ 3\mid 2\right] }(z;\alpha)
  =
 \det \begin{pmatrix}
    \cL_{0,1}& -1 & 0 & 0 & 0 \\[2pt]
    \cL_{0,3}& \cL_{1,2} & \cL_{2,1} & -1 & 0 \\[2pt]
    z^4\bcL_{0,3} & z^3\bcL_{1,3} &  z^2\bcL_{2,4} & z\bcL_{3,5} & \bcL_{4,6}\\[2pt]
    z^4\Lambda_{0,3} & z^3\Lambda_{1,3} & z^2\Lambda_{2,3}&
    z\Lambda_{3,3}&\Lambda_{4,3} \\[2pt]
    \bLambda_{0,2} & \bLambda_{1,2} & \bLambda_{2,2} & \bLambda_{3,2}
    &    \bLambda_{4,2}  
  \end{pmatrix}\,,\\
  L^{\left[ 3,2,1\mid \varnothing \right] \otimes \left[ 4,1^2\mid
      \varnothing\right] }(z;\alpha) 
  =
 \det \begin{pmatrix}
    \cL_{0,1}& -1 & 0 & 0 & 0 &0\\
    \cL_{0,3}& \cL_{1,2} & \cL_{2,1} & -1 & 0  &0\\
    \cL_{0,5}& \cL_{1,4} & \cL_{2,3} & \cL_{3,2} & \cL_{4,1}  & \cL_{5,0}\\
    z^5 \Lambda_{0,1} & z^4 \Lambda_{1,1} & z^3 \Lambda_{2,1} &  z^2\Lambda_{3,1} &
    z \Lambda_{4,1} & \Lambda_{5,1}  \\
    z^5 \Lambda_{0,2} & z^4 \Lambda_{1,2} & z^3 \Lambda_{2,2} &  z^2\Lambda_{3,2} &
    z \Lambda_{4,2} & \Lambda_{5,2}  \\
    z^5 \Lambda_{0,6} & z^4 \Lambda_{1,6} & z^3 \Lambda_{2,6} &  z^2\Lambda_{3,6} &
    z \Lambda_{4,6} & \Lambda_{5,6}  
  \end{pmatrix}\,.
\end{gather}

An explicit calculation applied to equation (\ref{EquivXLP}) now gives
\begin{equation}\label{eq:equivLex}
  \frac{1}{90} (\alpha-4)\alpha\ffac{\alpha+2}{4} \;L^{\left[ 2,1\mid 1\right]
    \otimes \left[ 3\mid 
      2\right] }\left( z;\alpha  
  \right) =  z^{-4}\; L^{\left[ 3,2,1\mid \varnothing \right] \otimes
    \lbrack 4,1^{2}\mid \varnothing ]}\left( z;\alpha+1\right) .
\end{equation}

\begin{figure}[ht]
\centering
\includegraphics{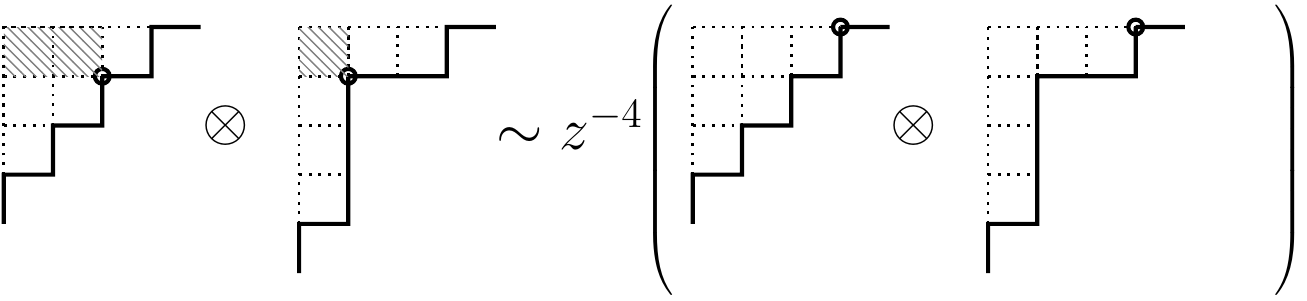}
\caption{Graphical description of the equivalence formula \eqref{eq:equivLex} between two Laguerre pseudo-Wronskians.}
\end{figure}

\end{example}

%%%%%%%%%%%%%%%%%%%%%%%%%%%%%%%%%%%%%%%%%%%%%%%%%%%%%%%%%%%%%%%%%%%
\section{Equivalence relations for the Jacobi pseudowronskians}\label{sec:J}
%%%%%%%%%%%%%%%%%%%%%%%%%%%%%%%%%%%%%%%%%%%%%%%%%%%%%%%%%%%%%%%%%%%

Proceeding exactly as in the case of the isotonic oscillator of the previous section, we obtain the following shifted Wronskians 
\begin{equation}
\left\{ 
\begin{array}{l}
W^{\left[ d\mid \overline{d}\right] \otimes \left[ \delta \mid \overline{%
\delta }\right] }\left( x;\alpha ,\beta \right) \propto \psi _{0}\left(
x;\alpha ,\beta \right) W^{_{1}\left[ d\mid \overline{d}\right] \otimes %
\left[ \delta \mid \overline{\delta }\right] }\left( x;\alpha _{1},\beta
_{1}\right) \\ 
W^{\left[ d\mid \overline{d}\right] \otimes \left[ \delta \mid \overline{%
\delta }\right] }\left( x;\alpha ,\beta \right) \propto \psi _{-1}\left(
x;\alpha ,\beta \right) W^{\left[ d\mid \overline{d}\right] _{1}\otimes %
\left[ \delta \mid \overline{\delta }\right] }\left( x;\alpha _{-1},\beta
_{-1}\right) \\ 
W^{\left[ d\mid \overline{d}\right] \otimes \left[ \delta \mid \overline{%
\delta }\right] }\left( x;\alpha ,\beta \right) \propto \phi _{0}\left(
x;\alpha ,\beta \right) W^{\left[ d\mid \overline{d}\right] \otimes _{1}%
\left[ \delta \mid \overline{\delta }\right] }\left( x;\alpha _{-1},\beta
_{1}\right) \\ 
W^{\left[ d\mid \overline{d}\right] \otimes \left[ \delta \mid \overline{%
\delta }\right] }\left( x;\alpha ,\beta \right) \propto \phi _{-1}\left(
x;\alpha ,\beta \right) W^{\left[ d\mid \overline{d}\right] \otimes \left[
\delta \mid \overline{\delta }\right] _{1}}\left( x;\alpha _{1},\beta
_{-1}\right)%
\end{array}%
\right.
\end{equation}%
and consequently

\begin{eqnarray}
W^{\left[ d\mid \overline{d}\right] \otimes \left[ \delta \mid \overline{%
\delta }\right] }\left( x;\alpha ,\beta \right) &\propto
&\prod_{i=0}^{d_{1}+k-1}\psi _{0}\left( x;\alpha _{i},\beta _{i}\right)
\prod_{j=0}^{\delta _{1}+\kappa -1}\phi _{0}\left( x;\alpha
_{-j+d_{1}+k},\beta _{j+d_{1}+k}\right) \\
&&\times W^{\left[ \varnothing \mid \overline{\lambda }\right] \otimes \left[
\varnothing \mid \overline{\mu }\right] }\left( x;\alpha _{d_{1}+k-(\delta
_{1}+\kappa )},\beta _{\delta _{1}+\kappa +d_{1}+k}\right) ,  \notag
\end{eqnarray}

\begin{eqnarray}
W^{\left[ d\mid \overline{d}\right] \otimes \left[ \delta \mid \overline{%
\delta }\right] }\left( x;\alpha ,\beta \right) &\propto &\prod_{i=0}^{\overline{%
d}_{1}+\overline{k}-1}\psi _{-1}\left( x;\alpha _{-i},\beta _{-i}\right)
\prod_{j=0}^{\overline{\delta }_{1}+\overline{\kappa }-1}\phi _{-1}\left(
x;\alpha _{j-(\overline{d}_{1}+\overline{k})},\beta _{-j-(\overline{d}_{1}+%
\overline{k})}\right) \\
&&\times W^{\left[ \lambda \mid \varnothing \right] \otimes \left[ \mu \mid
\varnothing \right] }\left( x;\alpha _{\overline{\delta }_{1}+\overline{%
\kappa }-(\overline{d}_{1}+\overline{k})},\beta _{-\left( \overline{d}_{1}+%
\overline{k}+\overline{\delta }_{1}+\overline{\kappa }\right) }\right) 
\notag
\end{eqnarray}%
and

\begin{eqnarray}\label{caca}
W^{\left[ d\mid \overline{d}\right] \otimes \left[ \delta \mid \overline{%
\delta }\right] }\left( x;\alpha ,\beta \right) &\propto &\prod_{i=0}^{\overline{%
d}_{1}+\overline{k}-1}\psi _{-1}\left( x;\alpha _{-i},\beta _{-i}\right)
\prod_{j=0}^{\overline{\delta }_{1}+\overline{\kappa }-1}\phi _{0}\left(
x;\alpha _{-j-(\overline{d}_{1}+\overline{k})},\beta _{j-(\overline{d}_{1}+%
\overline{k})}\right) \\
&&\times W^{\left[ \lambda \mid \varnothing \right] \otimes \left[
\varnothing \mid \overline{\mu }\right] }\left( x;\alpha _{-(\overline{d}%
_{1}+\overline{k}+\overline{\delta }_{1}+\overline{\kappa })},\beta _{-(%
\overline{d}_{1}+\overline{k})+\overline{\delta }_{1}+\overline{\kappa }%
}\right) .  \notag
\end{eqnarray}

For the extended potentials, we have

\begin{equation}
\left\{ 
\begin{array}{l}
V^{\left[ d\mid \overline{d}\right] \otimes \left[ \delta \mid \overline{%
\delta }\right] }\left( x;\alpha ,\beta \right) =V^{\left[ \lambda \mid
\varnothing \right] \otimes \left[ \mu \mid \varnothing \right] }\left(
x;\alpha _{d_{1}+k-(\delta _{1}+\kappa )},\beta _{d_{1}+k+\delta _{1}+\kappa
}\right) +E_{-\left( \overline{d}_{1}+\overline{k}\right) }(\alpha ,\beta )
\\ 
V^{\left[ d\mid \overline{d}\right] \otimes \left[ \delta \mid \overline{%
\delta }\right] }\left( x;\alpha ,\beta \right) =V^{\left[ \varnothing \mid 
\overline{\lambda }\right] \otimes \left[ \varnothing \mid \overline{\mu }%
\right] }(x;\alpha _{\delta _{1}+\kappa -(d_{1}+k)},\beta _{-(\delta
_{1}+\kappa +d_{1}+k)})+E_{d_{1}+k}(\alpha ,\beta ) \\ 
V^{\left[ d\mid \overline{d}\right] \otimes \left[ \delta \mid \overline{%
\delta }\right] }\left( x;\alpha ,\beta \right) =V^{\left[ \lambda \mid
\varnothing \right] \otimes \left[ \varnothing \mid \overline{\mu }\right]
}(x;\alpha _{d_{1}+k+\delta _{1}+\kappa },\beta _{d_{1}+k-(\delta
_{1}+\kappa )})+E_{-\left( \overline{d}_{1}+\overline{k}\right) }(\alpha
,\beta )\,.%
\end{array}%
\right.
\end{equation}

Using (\ref{conjshadtDPT}) and (\ref{conjtDPT}), expression (\ref{caca})
becomes

\begin{eqnarray*}
W^{\left[ d\mid \overline{d}\right] \otimes \left[ \delta \mid \overline{%
\delta }\right] }\left( x;\alpha,\beta \right) &\propto 
&(1-z)^{(\overline{d}_{1}+\overline{k}-\overline{\delta }_{1}-\overline{\kappa }%
)^{2}/4  -\alpha (\overline{d}_{1}+\overline{k}-%
\overline{\delta }_{1}-\overline{\kappa })/2} \cdot \left( 1+z\right) ^{(\overline{d}_{1}+\overline{k}+\overline{\delta }%
_{1}+\overline{\kappa })^{2}/4 -\beta (%
\overline{d}_{1}+\overline{k}+\overline{\delta }_{1}+\overline{\kappa }%
)/2} \\
&&\quad \times\, W^{\left[ \lambda \mid \varnothing \right] \otimes \left[ \mu \mid
\varnothing \right] }\left( x;\alpha _{-(\overline{d}_{1}+\overline{k})+%
\overline{\delta }_{1}+\overline{\kappa }},\beta _{-(\overline{d}_{1}+%
\overline{k}+\overline{\delta }_{1}+\overline{\kappa })}\right) \,,  \notag
\end{eqnarray*}%
where $z=\cos 2x$.

Using (\ref{wronskprop}) and the identities (\ref{conjtDPT}), (\ref
{shadtDPT}) and (\ref{conjshadtDPT}), we can write 
\begin{equation}
\left\{ 
\begin{array}{l}
\psi _{-\left( n+1\right) }\left( x;\alpha ,\beta \right) =\left( 1-z\right)
^{-\alpha }\left( 1+z\right) ^{-\beta }\psi _{0}\left( x;\alpha ,\beta
\right) \mathit{P}_{n}^{\left( -\alpha ,-\beta \right) }\left( z\right) \\ 
\phi _{n}\left( x;\alpha ,\beta \right) =\left( 1-z\right) ^{-\alpha }\psi
_{0}\left( x;\alpha ,\beta \right) \mathit{P}_{n}^{\left( -\alpha ,\beta
\right) }\left( z\right) \\ 
\phi _{-\left( n+1\right) }\left( x;\alpha ,\beta \right) =\left( 1+z\right)
^{-\beta }\psi _{0}\left( x;\alpha ,\beta \right) \mathit{P}_{n}^{\left(
\alpha ,-\beta \right) }\left( z\right)%
\end{array}%
\right. ,
\end{equation}%
which combined with the properties of the derivatives of Jacobi polynomials (\ref{derivJac102}), and noting that $m=k+\overline{k},r=\kappa +\overline{\kappa }$, gives

\begin{eqnarray*}
W^{\left[ d\mid \overline{d}\right] \otimes \left[ \delta \mid \overline{%
\delta }\right] }\left( x;\alpha ,\beta \right) &\propto &\left( 1-z\right)
^{\left( k-\overline{k}-\kappa +\overline{\kappa }\right) ^{2}/4-\left( 
\overline{k}+\kappa \right) \left( \overline{k}+\kappa -1\right) }\left(
1+z\right) ^{\left( k-\overline{k}+\kappa -\overline{\kappa }\right)
^{2}/4-\left( \overline{k}+\overline{\kappa }\right) \left( \overline{k}+%
\overline{\kappa }-1\right) }  \label{W1tDPT} \\
&&\times \left( 1-z\right) ^{\alpha \left( k-\overline{k}-\kappa +\overline{%
\kappa }\right) /2}\left( 1+z\right) ^{\beta \left( k-\overline{k}+\kappa -%
\overline{\kappa }\right) /2}\times P^{\left[ d\mid \overline{d}\right]
\otimes \left[ \delta \mid \overline{\delta }\right] }\left( z;\alpha ,\beta
\right) ,  \notag
\end{eqnarray*}%
where  $P^{\left[ d\mid \overline{d}\right] \otimes \left[ \delta \mid \overline{%
\delta }\right] }\left( z;\alpha ,\beta \right)$ is a polynomial which we shall denote as a \textit{Jacobi pseudowronskian}. A Jacobi pseudowronskian is defined by the following determinantal expression

\begin{equation}
P^{\left[ d\mid \overline{d}\right] \otimes \left[ \delta \mid \overline{%
\delta }\right] }\left( z;\alpha ,\beta \right) =\det \left( \left\{ 
\overrightarrow{\mathcal{P}}^{\left( j\right) } \right\} _{1\leq j\leq k},\left\{ \overrightarrow{\overline{\mathcal{%
P}}}^{\left( j\right) } \right\} _{1\leq j\leq 
\overline{k}},\left\{ \overrightarrow{\Pi }^{\left( j\right) } \right\} _{1\leq j\leq \kappa },\left\{ 
\overrightarrow{\overline{\Pi }}^{\left( j\right) }\right\} _{1\leq j\leq \overline{\kappa }}\right) ,  \label{JTtDPT}
\end{equation}%
where
\begin{equation}
\left\{ 
\begin{array}{l}
\mathcal{P}_{i}^{\left( j\right) }(z;\alpha ,\beta )=\frac{(n_{j}+\alpha
+\beta +1)_{i-1}}{2^{i-1}}\mathit{P}_{n_{j}-i+1}^{\left( \alpha +i-1,\beta
+i-1\right) }\left( z\right) \\ 
\overline{\mathcal{P}}_{i}^{\left( j\right) }(z;\alpha ,\beta )=\left(
-2\right) ^{i-1}\left( -\overline{n}_{j}\right) _{\overline{i-1}}\times
\left( 1-z^{2}\right) ^{m+r-i}\mathit{P}_{-\overline{n}_{j}+i-2}^{\left(
-\alpha -i+1,-\beta -i+1\right) }\left( z\right) \\ 
\Pi _{i}^{\left( j\right) }(z;\alpha ,\beta )=\left( -1\right) ^{i-1}\left(
l_{j}-\alpha \right) _{\underline{i-1}}\left( 1-z\right) ^{m+r-i}\emph{P}%
_{l_{j}}^{\left( -\alpha -i+1,\beta \right) }\left( z\right) \\ 
\overline{\Pi }_{i}^{\left( j\right) }(z;\alpha ,\beta )=\left( -\overline{l}%
_{j}-1-\beta \right) _{\underline{i-1}}\left( 1+z\right) ^{m+r-i}\emph{P}_{-%
\overline{l}_{j}-1}^{\left( \alpha ,-\beta -i+1\right) }\left( z\right) ,%
\end{array}%
\right.
\end{equation}%
with $n_{j}=d_{j}+k-j\geq 0$, $l_{j}=\delta _{j}+\kappa -j\geq 0$, $-%
\overline{n}_{j}-1=\overline{d}_{j}+\overline{k}-j\geq 0$ and $-\overline{l}%
_{j}-1=\overline{\delta }_{j}+\overline{\kappa }-j\geq 0$.

If $\left[ \lambda \mid \varnothing \right] $ is in the same shape class as $%
\left[ d\mid \overline{d}\right] =\left[ d_{1},...,d_{k}\mid \overline{d}%
_{1},...,\overline{d}_{\overline{k}}\right] $, we can use (\ref{W1tDPT}) replacing $\overline{k}$ and $\overline{\kappa }$ by $0$ and $k$ and $%
\kappa $ by $l(\lambda )=\overline{d}_{1}+k$ and $l(\mu )$ $=\overline{%
\delta }_{1}+\kappa $ respectively. This gives

\begin{eqnarray}
W^{\left[ \lambda \mid \varnothing \right] \otimes \left[ \mu \mid
\varnothing \right] }\left( x;\alpha ,\beta \right) &\propto &\left( 1-z\right)
^{\left( \overline{d}_{1}+k-\overline{\delta }_{1}-\kappa \right)
^{2}/4-\left( \overline{\delta }_{1}+\kappa \right) \left( \overline{\delta }%
_{1}+\kappa -1\right) }\left( 1+z\right) ^{\left( \overline{d}_{1}+k+%
\overline{\delta }_{1}+\kappa \right) ^{2}/4}  \label{W1tDPT12} \\ \notag
&&\quad \times \left( 1-z\right) ^{\alpha \left( \overline{d}_{1}+k-\overline{%
\delta }_{1}-\kappa \right) /2}\left( 1+z\right) ^{\beta \left( \overline{d}%
_{1}+k+\overline{\delta }_{1}+\kappa \right) /2}\\
&& \quad  \times \, P^{\left[ \lambda
\mid \varnothing \right] \otimes \left[ \mu \mid \varnothing \right] }\left(
z;\alpha ,\beta \right) .  \notag
\end{eqnarray}
The equivalence relation can then be rewritten as

\begin{eqnarray}
W^{\left[ d\mid \overline{d}\right] \otimes \left[ \delta \mid \overline{%
\delta }\right] }\left( x;\alpha ,\beta \right) &\propto &\left( 1-z\right)
^{(k-\overline{k}-\kappa +\overline{\kappa })^{2}/4-\left( \overline{\delta }%
_{1}+\kappa \right) \left( \overline{\delta }_{1}+\kappa -1\right) }\left(
1+z\right) ^{(k-\overline{k}+\kappa -\overline{\kappa })^{2}/4} \\ \notag
&&\quad \times \left( 1-z\right) ^{\alpha (k-\overline{k}-\kappa +\overline{\kappa 
})/2}\left( 1+z\right) ^{\beta (k-\overline{k}+\kappa -\overline{\kappa }%
)/2} \\ \notag
&& \quad \times P^{\left[ \lambda \mid \varnothing \right] \otimes \left[ \mu
\mid \varnothing \right] }\left( z;\alpha _{-(\overline{d}_{1}+\overline{k})+%
\overline{\delta }_{1}+\overline{\kappa }},\beta _{-(\overline{d}_{1}+%
\overline{k}+\overline{\delta }_{1}+\overline{\kappa })}\right)
\end{eqnarray}%
By making use of (\ref{W1tDPT}), the previous formula can be finally expressed as

\begin{eqnarray}
P^{\left[ d\mid \overline{d}\right] \otimes \left[ \delta \mid \overline{%
\delta }\right] }\left( z;\alpha ,\beta \right)  &\propto &\left( 1-z\right)
^{\left( \overline{k}+\kappa \right) \left( \overline{k}+\kappa -1\right)
-\left( \overline{\delta }_{1}+\kappa \right) \left( \overline{\delta }%
_{1}+\kappa -1\right) }\left( 1+z\right) ^{\left( \overline{k}+\overline{%
\kappa }\right) \left( \overline{k}+\overline{\kappa }-1\right) }
\label{EquivXJP} \\
&&\times P^{\left[ \lambda \mid \varnothing \right] \otimes \left[ \mu \mid
\varnothing \right] }\left( z;\alpha _{-(\overline{d}_{1}+\overline{k})+%
\overline{\delta }_{1}+\overline{\kappa }},\beta _{-(\overline{d}_{1}+%
\overline{k}+\overline{\delta }_{1}+\overline{\kappa })}\right) ,  \notag
\end{eqnarray}%
which is the \textit{equivalence relation between Jacobi pseudowronskians}.

\begin{example}

As an example of an equivalence relation, consider the case where $\left[ d\mid 
\overline{d}\right] =\left[ 2,1\mid 1\right] $, $N_{3}=(3,1,-2)$ and $\left[
\delta \mid \overline{\delta }\right] =\left[ 3\mid 2\right] $, $%
L_{2}=(3,-3) $. We then have $\lambda =\left( 3,2,1\right) $ and $\mu
=\left( 4,1^{2}\right) $. 

Setting 
\begin{align}
  \cP_{i,j} &= 2^{-i}\rfac{j+\alpha+\beta+1}{i} P^{(\alpha+i,\beta+i)}_{j-i} \\
  \bcP_{i,j} &= 2^i \rfac{j+1}{i} P^{(-\alpha-i,-\beta-i)}_{j+i} \\
  \Pi_{i,j} &=   \rfac{j-\alpha}{i} P^{(-\alpha-i,\beta+i)}_j\\
  \bPi_{i,j} &=  \rfac{j-\beta}{i} P^{(\alpha+i,-\beta-i)}_j
\end{align}
 the pseudo-Wronskians in question have the explicit
determinantal form
\begin{gather}
  P^{\left[ 2,1\mid 1\right] \otimes \left[ 3\mid 2\right] }(z;\alpha,\beta)
  =
 \det \begin{pmatrix}
    \cP_{0,1}& \cP_{1,0} & 0 & 0 & 0 \\[2pt]
    \cP_{0,3}& \cP_{1,2} & \cP_{2,1} & \cP_{3,0} & 0 \\[2pt]
    (1-z^2)^4\bcP_{0,3} & (1-z^2)^3\bcP_{1,3} &  (1-z^2)^2\bcP_{2,4} &
    (1-z^2)\bcP_{3,5} & \bcP_{4,6}\\[2pt]  
    (1-z)^4\Pi_{0,3} & (1-z)^3\Pi_{1,3} & (1-z)^2\Pi_{2,3}&
    (1-z)\Pi_{3,3}&\Pi_{4,3} \\[2pt]
    (1+z)^4\bPi_{0,2} & (1+z)^3\bPi_{1,2} & (1+z)^2\bPi_{2,2} & (1+z)\bPi_{3,2}
    &    \bPi_{4,2}  
  \end{pmatrix} \,,\\
  P^{\left[ 3,2,1\mid \varnothing \right] \otimes \left[ 4,1^2\mid
      \varnothing\right] }(z;\alpha,\beta) 
  =
 \det \begin{pmatrix}
    \cP_{0,1}& \cP_{1,0} & 0 & 0 & 0 &0\\
    \cP_{0,3}& \cP_{1,2} & \cP_{2,1} & \cP_{3,0} & 0  &0\\
    \cP_{0,5}& \cP_{1,4} & \cP_{2,3} & \cP_{3,2} & \cP_{4,1}  & \cP_{5,0}\\
    (1-z)^5 \Pi_{0,1} & (1-z)^4 \Pi_{1,1} & (1-z)^3 \Pi_{2,1} &  (1-z)^2\Pi_{3,1} &
    (1-z) \Pi_{4,1} & \Pi_{5,1}  \\
    (1-z)^5 \Pi_{0,2} & (1-z)^4 \Pi_{1,2} & (1-z)^3 \Pi_{2,2} &  (1-z)^2\Pi_{3,2} &
    (1-z) \Pi_{4,2} & \Pi_{5,2}  \\
    (1-z)^5 \Pi_{0,6} & (1-z)^4 \Pi_{1,6} & (1-z)^3 \Pi_{2,6} &  (1-z)^2\Pi_{3,6} &
    (1-z) \Pi_{4,6} & \Pi_{5,6}  
  \end{pmatrix}\,.
\end{gather}
An explicit calculation applied to equation (\ref{EquivXJP}) now gives
\begin{equation}\label{eq:equivJex}
K(\alpha,\beta) P^{\left[ 2,1\mid 1\right] \otimes \left[ 3\mid 2\right] }\left( z;\alpha
,\beta \right) = \left( 1-z\right) ^{-4}\left( 1+z\right) ^{2}\times P^{%
\left[ 3,2,1\mid \varnothing \right] \otimes \lbrack 4,1^{2}\mid \varnothing
]}\left( z;\alpha _{1},\beta _{-5}\right) ,
\end{equation}%
where
\[ K(\alpha,\beta)= \frac{23040 (\beta-1)(\beta+1)(\beta+2)}{(\alpha-4)\alpha
  \ffac{\alpha+2}{4}\ffac{\alpha-\beta-3}{2} (\alpha-\beta+2)(\alpha+\beta+1)(\alpha+\beta+3)}.\]

\begin{figure}[ht]
\centering
\includegraphics{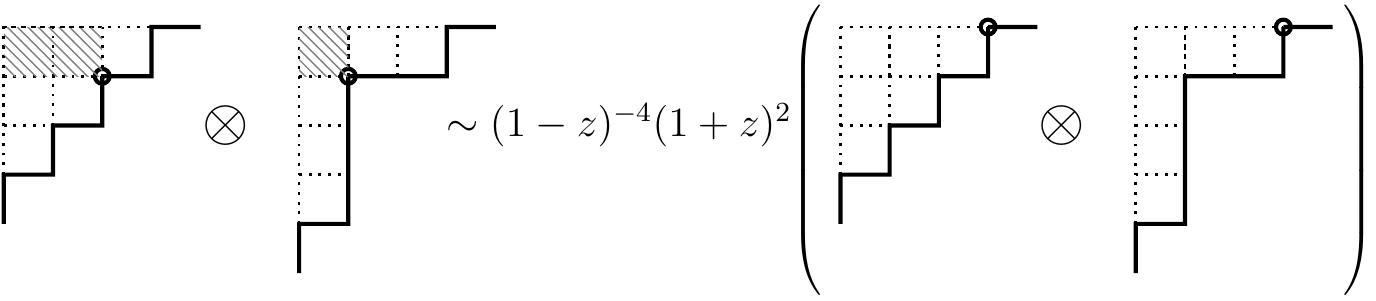}
\caption{Graphical description of the equivalence formula \eqref{eq:equivJex} between two Jacobi pseudo-Wronskians.}
\end{figure}

% \todo{ Illustrate this identity graphically replacing the P
% Wronskians with punctured young diagrams (joined by $\otimes$). Maintain the 
% $(1-z)^{-4}(1+z)^2$ symbols}

\end{example}

%%%%%%%%%%%%%%%%%%%%%%%%%%%%%%%%%%%%%%%%%%%%%%%%%%%%%%%%%%%%%%%%%%%
\section{Summary and Outlook}
%%%%%%%%%%%%%%%%%%%%%%%%%%%%%%%%%%%%%%%%%%%%%%%%%%%%%%%%%%%%%%%%%%%

In \cite{GGM4} we introduced the concept of Hermite pseudo-Wronskians determinants and proved many equivalence formulas among them. In this paper we have extended the analysis to Laguerre and Jacobi pseudo-Wronskians.
Equivalence formulas express different but equivalent manners to perform rational Darboux transformations on the isotonic oscillator and tDPT potential such that the final potential is the same up to an additive shift.

In the Laguerre and Jacobi setting, the situation is richer, as rational Darboux transformations can be state-adding, state-deleting or isospectral, i.e. we can choose seed functions for rational Darboux transformations out of four families, instead of just two families in the Hermite case. Thus, equivalences are obtained by shifting the origin in two independent Maya diagrams.   The equivalence relations  can be viewed as a
generalization of the usual shape invariance property at a multi-step level, combined with the discrete $\Gamma$ symmetries of the primary potentials (isotonic and tDPT).

As discussed in \cite{GGM4} for the generalized Hermite and Okamoto classes, the existence of these equivalence relations provides optimized determinantal representations for the polynomials
associated to rational solutions of (higher order) Painlev\'{e} equations of $A_N$-type. This will be the subject of further investigation, aimed at providing a complete classification of the rational solutions of $A_N$-Painlev\'e equations, together with their optimized (pseudo)-Wronskian representation.

\section*{Acknowledgements}

The research of D.G.U. has been supported in
part by Spanish MINECO-FEDER Grant 
MTM2015-65888-C4-3. He also acknowledges financial support from the Spanish Ministry of Economy and Competitiveness, through the ``Severo Ochoa Programme for Centres of Excellence in R\&D'' (SEV-2015-0554).
The research of R.M. was supported in
part by NSERC grant RGPIN-228057-2009. They both would like to thank
 Universit\'e de Lorraine for its hospitality during their visit in the summer of 2016 where this work was initiated. The draft was finalized during the thematic trimester on ``Orthogonal polynomials and Special functions  in approximation theory and mathematical physics'', held in Madrid at the Institute of Mathematical Sciences in the Fall of 2017, whose financial support is also acknolwedged.

\end{document}